\documentclass[10pt,a4paper]{amsart}

\usepackage{amsmath, amssymb, amsthm}
\usepackage{geometry}
\newtheorem{theorem}{Theorem}[section]
\newtheorem{lemma}[theorem]{Lemma}
\newtheorem{proposition}[theorem]{Proposition}

\theoremstyle{definition}
\newtheorem{definition}[theorem]{Definition}
\newtheorem{hypothesis}[theorem]{Hypothesis}

\theoremstyle{remark}
\newtheorem{remark}[theorem]{Remark}

\title[Rate-Preserving Perturbations for Sublinear ODEs]{Characterisation of Rate-Preserving Perturbations of First-Order Sublinear Differential Equations}

\author{John A. D. Appleby}
\address{School of Mathematical Sciences, Dublin City University, Glasnevin, Dublin 9, Ireland}
\email{john.appleby@dcu.ie}

\author{Subham Pal}
\address{School of Engineering, Trinity College Dublin, Dublin 2, Ireland}
\email{palsu@tcd.ie}
\address{School of Mathematical Sciences, Dublin City University, Glasnevin, Dublin 9, Ireland}
\email{subham.pal2@mail.dcu.ie}


\date{}

\begin{document}
	
	\begin{abstract}
		We study positive solutions of $x'(t)=-f(x(t))+g(t)$, where $g$ is a
		nonnegative forcing and $f$ is sublinear at zero. The unperturbed equation has
		a non-hyperbolic stable equilibrium, and its decay is described by the
		escape-time transform $F(x)=\int_x^1 f(u)^{-1}\,du$. We characterise
		preservation of this decay rate through the condition $F(x(t))/t\to1$.
		
		When the transformed nonlinearity $\varphi=f\circ F^{-1}$ is weakly
		asymptotic preserving (WAP), we prove necessary and sufficient conditions on
		$g$: a Ces\`aro condition, an equivalent finite-window integral condition,
		and, under a two-sided elasticity hypothesis, an equivalent tail-integral
		condition. We also characterise simultaneous preservation of the asymptotic
		behaviour of the solution and its derivative. The proofs use two
		Ces\`aro--Tauberian equivalences for monotone functions. Finally, we give
		checkable sufficient conditions for WAP in terms of
		$Q=ff''/(f')^2$ and construct examples showing how the characterisation
		changes when the structural hypotheses fail.
	\end{abstract}
	
	\maketitle
	
	\section{Introduction}
	\label{sec:introduction}
	
	We study the preservation of the asymptotic decay rate of positive solutions to
	the non-autonomous scalar ordinary differential equation
	\begin{equation}\label{eq:main_ode}
		x'(t)=-f(x(t))+g(t),\qquad t\geq0,\qquad x(0)=\xi>0,
	\end{equation}
	where $g$ is a continuous nonnegative forcing term and the unperturbed
	autonomous equation
	\[
	u'(t)=-f(u(t))
	\]
	has a non-hyperbolic asymptotically stable equilibrium at the origin.
	Throughout the main results, we work on the positive half-line and assume that
	$f(0)=0$, $f(x)>0$ for $x>0$, $f$ is increasing, and $f$ is sublinear at zero
	in the sense that $f(x)/x \to 0$ as $x \to 0^+$. Under these assumptions, the
	escape-time (or nonlinear clock) function
	\begin{equation}\label{eq:F}
		F(x)=\int_x^1\frac{du}{f(u)}
	\end{equation}
	satisfies $F(x)\to\infty$ as $x\to0^+$. Hence $F$ is strictly decreasing and
	admits a strictly decreasing inverse on its range, which contains a half-line
	$[T,\infty)$ for some $T>0$. Throughout, expressions involving $F^{-1}(t)$
	are understood for sufficiently large $t$; its definition or extension on a
	bounded initial interval has no effect on any asymptotic statement.
	
	For the unperturbed equation with initial value $u(0)=\zeta>0$, separation of
	variables gives $F(u(t))=F(\zeta)+t$ for all $t\geq0$, and therefore
	\[
	\lim_{t\to\infty}\frac{F(u(t))}{t}=1,
	\]
	independently of $\zeta>0$. We adopt this limit as a generalised nonlinear
	Lyapunov exponent and say that a perturbed solution $x$ of \eqref{eq:main_ode}
	preserves the unperturbed decay rate when
	\[
	\lim_{t\to\infty}\frac{F(x(t))}{t}=1.
	\]
	This criterion remains informative precisely in the sublinear, non-hyperbolic
	setting where solutions decay at power-law or slower-than-power-law rates and
	the classical exponential Lyapunov exponent degenerates to zero. Moreover, as
	explained in Section~\ref{sec:main_results}, the escape-time criterion
	$F(x(t))/t \to 1$ enjoys robustness properties---including initial-condition
	independence and forced integrability of $g$---that the seemingly more direct
	asymptotic equivalence $x(t)\sim u(t)$ can fail to possess for super-flat
	nonlinearities.
	
	\subsection{Principal results and contributions}
	
	The main result is a necessary-and-sufficient characterisation of rate
	preservation when $f$ is sublinear at zero and
	\[
	\varphi(t):=(f\circ F^{-1})(t)
	\]
	is weakly asymptotic preserving (WAP) at infinity; see
	Definition~\ref{def:WAP}. Theorem~\ref{thm.22} proves the equivalence of
	\begin{enumerate}
		\item[(A)] $F(x(t))/t\to1$;
		\item[(B)] $g\in L^1(0,\infty)$ and
		$\displaystyle \frac1t\int_0^t\frac{g(s)}{\varphi(s)}\,ds\to0$;
		\item[(C)] $g\in L^1(0,\infty)$ and
		$\displaystyle \frac{1}{t\varphi(t)}\int_t^{2t}g(s)\,ds\to0$.
	\end{enumerate}
	Under the additional two-sided elasticity condition~\eqref{eq.8},
	Theorem~\ref{thm.23} shows that these conditions are also equivalent to
	$\int_t^\infty g(s)\,ds=o(F^{-1}(t))$. Theorems~\ref{thm.24} and
	\ref{thm.24bis} identify the stronger pointwise condition on $g$ that
	characterises simultaneous preservation of the asymptotic behaviour of the
	solution and its derivative.
	
	The proof of the main characterisation requires two real-variable results.
	Theorems~\ref{thm.generalcesaro} and~\ref{thm.converselhopital} show, under
	WAP and a two-sided elasticity bound where appropriate, that a Ces\`aro
	condition on $a'/b'$ is equivalent first to a dyadic increment condition on
	$a$ and then to $a=o(b)$. The latter is a Ces\`aro-averaged converse to
	l'H\^opital's rule; it does not assert the generally false pointwise converse.
	
	A second issue is the verifiability of the hypothesis on the transformed
	function $\varphi$. We derive sufficient conditions directly in terms of
	$f$ and its derivatives, using
	\[
	Q(x)=\frac{f(x)f''(x)}{(f'(x))^2}.
	\]
	In particular, two-sided positive finite bounds on $Q$ imply WAP, and the
	criteria cover regularly varying and rapidly varying nonlinearities as well
	as useful monotonicity and log-concavity conditions.
	
	Finally, the examples in Section~\ref{sec:sharpness_counterexamples} and
	Appendices~A--D identify the limits of the characterisation. They show how
	the forcing criterion changes when WAP or sublinearity fails, and why the
	tail-integral condition requires the additional elasticity hypothesis. The
	results build on, rather than replace, the rate-preservation theory for
	regularly varying nonlinearities developed by Appleby and Patterson; the
	specific advance here is the role of monotone O-regular variation of
	$f\circ F^{-1}$, together with the associated Tauberian equivalences,
	verifiable derivative criteria, and boundary examples.
	
	\subsection{Organisation of the paper}
	
	Section~\ref{sec:lit_review} discusses the relevant literature. Section~\ref{sec:main_results} introduces the escape-time scale and WAP hypothesis,
	and states the main theorems. Section~\ref{sec:sharpness_counterexamples}
	summarises the boundary examples, whose full constructions appear in
	Appendices~A--D. Section~\ref{sec:cesaro_char} establishes the
	Ces\`aro--Tauberian equivalences. Section~\ref{sec:equiv_AB} proves the coupled ODE
	characterisation and applies the real-variable results to complete the proofs
	of Theorems~\ref{thm.22} and~\ref{thm.23}. Section~\ref{sec:sufficient_WAP}
	gives derivative criteria for WAP, and Section~\ref{sec:properties_WAP}
	collects auxiliary properties of WAP functions.
	
	\section{Literature and relation to previous work}
	\label{sec:lit_review}
	
	The escape-time transform $F(x)=\int_x^1 f(u)^{-1}\,du$ and its inverse are
	classical instruments in the asymptotic analysis of scalar differential
	equations. For hyperbolic equilibria, classical perturbation theory---as in
	Hartman~\cite{hartman2002ordinary} and Coppel~\cite{coppel1965stability}---relies
	on linearisation and variation-of-constants estimates to preserve exponential
	decay rates. Near non-hyperbolic equilibria, however, linearisation fails and
	regular variation has long provided the natural framework for deriving exact
	decay rates; see Avakumovi\'c~\cite{avakumovic1947sur}, Mari\'c and
	Tomi\'c~\cite{maric1976asymptotic,maric1980asymptotics},
	Omey~\cite{omey1981regular}, and the monographs of Bingham, Goldie, and
	Teugels~\cite{bingham1989regular} and Mari\'c~\cite{maric2000regular}.
	Subsequent developments have applied regular and rapid variation to sublinear
	Emden--Fowler equations, functional differential equations, cyclic systems, and
	higher-order non-autonomous equations; see
	\cite{kusano2006class,kusano2011precise,evtukhov2011asymptotic,
		kozma2012asymptotic,matucci2014asymptotics,matucci2018strongly,
		evtukhov2019asymptotics,chepok2020asymptotic}.
	
	The closest antecedents to the present study are the works of Appleby and
	Patterson on additive perturbations of scalar equations. In
	\cite{appleby2014necessary}, necessary and sufficient conditions for preserving
	decay rates were established for regularly varying nonlinearities, including
	the principal sublinear regime $f\in\mathrm{RV}_0(\beta)$ with $\beta>1$,
	alongside characterisations of derivative preservation and stochastic
	analogues. Their subsequent classification paper
	\cite{appleby2016classification} analysed how perturbations of different
	asymptotic sizes modify the decay rate and already identified exceptional
	behaviour in the borderline class $f\in\mathrm{RV}_0(1)$, where pointwise
	smallness of $g/(f\circ F^{-1})$ need not be necessary for rate preservation.
	In \cite{appleby2017exact}, related rate-preservation criteria were extended to
	equations with unbounded delay and asymptotic preserving nonlinearities. The
	present paper builds directly on these precursors rather than claiming that
	rate preservation for regularly varying nonlinearities is new.
	
	The present paper extends this line of work in three related directions.
	First, it moves beyond regular variation of $f$ to the requirement that
	$\varphi = f\circ F^{-1}$ belong to the monotone subclass of O-regularly
	varying functions ($OR$). The $OR$ class was introduced by
	Avakumovi\'c and Karamata and developed systematically by Aljan\v{c}i\'c and
	Arandelovi\'c~\cite{aljancic1977oregular} and Bingham, Goldie, and
	Teugels~\cite[Chapter~2]{bingham1989regular}, where growth bounds are governed
	by the Matuszewska indices (indeed, Lemma~\ref{lemma.wapmeansnofasterthanpowerlaw}
	reflects the finiteness of the lower Matuszewska index). Our contribution is to
	identify monotone O-regular variation of $f\circ F^{-1}$ as the exact analytic
	property required by the nonlinear perturbation argument, and to supply
	checkable differential criteria for it via $Q(x)=f(x)f''(x)/(f'(x))^2$---a
	quantity that also plays a central role in de~Haan's $\Gamma$-class and von
	Mises conditions for rapid variation~\cite[Chapter~3]{bingham1989regular}.
	
	Second, to decouple the forcing $g$ from $\varphi=f\circ F^{-1}$ inside the
	Ces\`aro condition of Theorem~\ref{thm.22}, we establish Tauberian-type
	equivalences for monotone functions (Section~\ref{sec:cesaro_char}). Real-variable
	Tauberian theorems linking Ces\`aro averages, dyadic increments, and ratios of
	monotone functions form a rich classical subject (see
	Korevaar~\cite{korevaar2004tauberian} and \cite{bingham1989regular}); our
	Theorems~\ref{thm.generalcesaro} and~\ref{thm.converselhopital} isolate the
	exact O-regular variation and two-sided elasticity hypotheses under which
	Ces\`aro convergence of $a'/b'$ is equivalent to $a(t)-a(2t)=o(t|b'(t)|)$ and
	to $a(t)=o(b(t))$.
	
	Third, the forcing enters additively and may have irregular local excursions
	even when its accumulated mass is small. This differs from equations with
	multiplicative time-dependent coefficients, studied for example by Evtukhov
	and collaborators~\cite{evtukhov2011asymptotic,evtukhov2019asymptotics}.
	The finite-window and tail-integral criteria capture this cumulative forcing
	mechanism; Appendices~A--D describe the limits of those criteria.
	
	\section{Rate scale, WAP hypothesis, and main results}
	\label{sec:main_results}
	
	In this section we examine the asymptotic behaviour of positive solutions of
	the ordinary differential equation \eqref{eq:main_ode}. We work with positive
	solutions and assume that $g(t)\geq0$ for all $t\geq0$ and $x(0)=\xi>0$. On
	the positive half-line we assume $f(0)=0$, $f(x)>0$ for $x>0$, and that $f$ is
	increasing. These hypotheses ensure that the unperturbed equation has an
	asymptotically stable equilibrium at zero on the positive half-line; the
	results are local near zero and do not require a whole-line formulation of
	global stability.
	
	When $g(t)\equiv 0$ for all $t\geq 0$, the equation is separable. Let $u$
	denote the unperturbed solution with initial condition $u(0)=\zeta>0$. Separation
	of variables gives
	\[
	F(u(t))=F(\zeta)+t, \qquad t\geq 0,
	\]
	where $F$ is the escape-time function \eqref{eq:F}:
	\[
	F(x)=\int_x^1 \frac{1}{f(u)}\,du.
	\]
	We are interested in the non-hyperbolic sublinear regime $f(x)/x\to 0$ as
	$x\to 0^+$ (for instance, $f'(0^+)=0$), in which case
	\[
	F(x)\to \infty \qquad \text{as } x\to 0^+,
	\]
	and the maximal interval of existence of $u$ is $[0,\infty)$. The function $F$
	is strictly decreasing and therefore has a strictly decreasing inverse $F^{-1}$
	obeying $F^{-1}(t)\to 0$ as $t\to\infty$. In the context of the dynamical
	system, it is natural to assume that $f$ is increasing near zero; for
	simplicity, we assume throughout the main theorems that $f$ is increasing on
	$[0,\infty)$.
	
	\subsection{Justification of the escape-time rate criterion}
	
	Although one can write $u(t)=F^{-1}(t+F(\zeta))$, it is desirable to measure
	decay-rate preservation in a way that is robust both to small perturbations and
	to changes in initial conditions. Tracking the limit of $F(x(t))/t$ as
	$t\to\infty$ achieves this. For example, suppose $z$ obeys the state-perturbed
	equation
	\[
	z'(t)=-f(z(t))+\phi(z(t)),\qquad t\geq 0,
	\]
	where $\phi(x)=o(f(x))$ as $x\to 0^+$. In particular, suppose $|\phi(x)|<f(x)$
	for all $x\in (0,\delta)$. If $z(0)=\xi\in (0,\delta)$, then $z(t)\to 0$ as
	$t\to\infty$ and $z(t)>0$ for all $t\geq 0$, so
	\[
	\lim_{t\to\infty} \frac{z'(t)}{f(z(t))}=-1.
	\]
	Hence, for every $\epsilon\in (0,1)$, there exists $T(\epsilon)>0$ such that
	\[
	-1-\epsilon< \frac{z'(t)}{f(z(t))}<-1+\epsilon, \qquad t\geq T(\epsilon).
	\]
	Integrating over $[T,t]$ and noting that
	$\int_{T}^t \frac{z'(s)}{f(z(s))}\,ds=-F(z(t))+F(z(T))$ gives
	\[
	(1-\epsilon)(t-T)+F(z(T))\leq F(z(t))\leq (1+\epsilon)(t-T)+F(z(T)), \qquad t\geq T.
	\]
	Dividing by $t$, letting $t\to\infty$, and then letting $\epsilon\to 0^+$ yields
	\[
	\lim_{t\to\infty} \frac{F(z(t))}{t}=1.
	\]
	When $f(x)=ax$ and $\phi(x)=o(x)$, this reduces to the statement that the
	nonlinear equation has the same Lyapunov exponent as $u'=-au$. Thus the limit
	of $F(x(t))/t$ is a natural nonlinear generalisation of the Lyapunov exponent,
	and is automatically independent of the initial condition because
	$\lim_{t\to\infty} F(u(t))/t=1$ for every $\zeta>0$.
	
	A seemingly more direct way to measure rate preservation is to require
	$x(t)\sim u(t)$ as $t\to\infty$. For superlinear equations ($f(x)/x\to\infty$
	as $x\to 0^+$), however, even unperturbed solutions fail to share a common
	asymptotic profile: if $\xi_1>\xi_2$, then $u(t;\xi_1)/u(t;\xi_2)\to\infty$ as
	$t\to\infty$. In the sublinear case, all unperturbed solutions do satisfy
	$u(t;\xi)\sim u(t;1)$, yet the direct equivalence $x(t)\sim u(t)$ suffers from
	a different drawback: as our main theorems show,
	$F(x(t))/t\to 1$ always forces $g\in L^1(0,\infty)$ (so the cumulative forcing
	acts like a bounded shift of the initial condition), whereas $x(t)\sim F^{-1}(t)$
	need \emph{not} imply $g\in L^1(0,\infty)$ when $f$ is super-flat and $F^{-1}$
	is slowly varying. An explicit example of this phenomenon with
	$f(x)=x^2 e^{-1/x}$ and $\varphi = f\circ F^{-1}\in \mathrm{RV}_\infty(-1)$
	is established in Proposition~\ref{prop:solution-asymptotic-not-tail}
	(Appendix~\ref{app:superflat}). By contrast, the criterion $F(x(t))/t\to 1$
	behaves uniformly across all classes of nonlinearities considered here: it is
	invariant under changes of initial condition, stable under subdominant
	state-dependent perturbations, and compatible with integrability of additive
	forcing.
	
	\subsection{Weakly asymptotic preserving functions}
	
	Granted that $F(x(t))/t\to 1$ is the appropriate measure of rate preservation,
	we seek necessary and sufficient smallness conditions on $g$ ensuring that $x$
	inherits this rate. Roughly speaking, our main results cover the entire
	spectrum of weak nonlinearities between exceptionally flat functions
	($f^{(n)}(0)=0$ for all $n\in\mathbb{N}$) and power-law nonlinearities
	$f(x)\sim x^\beta$ ($\beta>1$), stopping precisely at the almost-linear
	boundary $f\in\mathrm{RV}_0(1)$ where unperturbed solutions decay faster than
	any power law.
	
	To formulate the governing structural hypothesis on $f$, we define the
	following class of functions.
	
	\begin{definition}[Weakly Asymptotic Preserving Function] \label{def:WAP}
		Let $\varphi : (0,\infty) \to (0,\infty)$ be a function.
		\begin{enumerate}
			\item[(i)] $\varphi$ is \emph{weakly asymptotic preserving (WAP) at $0$} if it is non-decreasing on $(0,\delta)$ for some $\delta > 0$ and if
			\begin{equation} \label{eq.WAPlimits0}
				0 < \underline{\Phi}_\lambda :=
				\liminf_{x \to 0^+} \frac{\varphi(\lambda x)}{\varphi(x)}
				\le
				\limsup_{x \to 0^+} \frac{\varphi(\lambda x)}{\varphi(x)}
				=: \overline{\Phi}_\lambda < +\infty,
			\end{equation}
			for all $\lambda \in (0,\infty)$.
			
			\item[(ii)] $\varphi$ is \emph{weakly asymptotic preserving (WAP) at $\infty$} if it is non-increasing on $(T,\infty)$ for some $T > 0$ and if
			\begin{equation} \label{eq.WAPlimitsInf}
				0 < \underline{\Phi}_\lambda
				:=\liminf_{x \to \infty} \frac{\varphi(\lambda x)}{\varphi(x)}
				\le
				\limsup_{x \to \infty} \frac{\varphi(\lambda x)}{\varphi(x)}
				=: \overline{\Phi}_\lambda < +\infty,
			\end{equation}
			for all $\lambda \in (0,\infty)$.
		\end{enumerate}
	\end{definition}
	
	For positive functions, the dilation-ratio condition \eqref{eq.WAPlimitsInf} is
	the defining O-regular variation condition at infinity (class $OR$; see
	\cite[Chapter~2]{bingham1989regular}), combined with the eventual monotonicity
	used in our ODE estimates. Thus WAP at $\infty$ is the non-increasing subclass
	of $OR$. We retain the terminology WAP to highlight its role in preserving
	asymptotic decay rates under additive perturbations.
	
	We record two auxiliary properties of WAP functions, proved in
	Section~\ref{sec:properties_WAP}.
	
	\begin{lemma}\label{lemma.intwapiswap}
		Suppose $\varphi$ is WAP at $\infty$ and $\varphi \in L^1(0,\infty)$. Then
		\[
		\phi(t) := \int_t^\infty \varphi(s)\,ds
		\]
		is also WAP at $\infty$.
	\end{lemma}
	
	\begin{lemma} \label{lemma.wapmeansnofasterthanpowerlaw}
		Let $\varphi: (0, \infty) \to (0, \infty)$ be a decreasing function that is WAP at $\infty$. Define the limit inferior of its logarithmic decay rate by
		\[
		\lambda = \liminf_{t\to\infty} \frac{\log \varphi(t)}{\log t}.
		\]
		Then $\lambda > -\infty$.
	\end{lemma}
	
	The standing structural hypothesis of our main theorems is that
	\emph{$f\circ F^{-1}$ is WAP at $\infty$}. This is not automatic from
	monotonicity alone, since monotonicity provides only one-sided ratio bounds.
	Note that $\varphi = f\circ F^{-1}$ is automatically integrable near $\infty$,
	because the substitution $u = F^{-1}(s)$ yields
	\[
	\int_t^\infty \varphi(s)\,ds = F^{-1}(t).
	\]
	Consequently, whenever $f\circ F^{-1}$ is WAP at $\infty$,
	Lemma~\ref{lemma.intwapiswap} implies that $F^{-1}$ is also WAP at $\infty$,
	and Lemma~\ref{lemma.wapmeansnofasterthanpowerlaw} applied to $F^{-1}$ shows
	that unperturbed solutions decay no faster than a power law.
	
	Many natural sublinear nonlinearities satisfy this WAP condition. Every
	monotone regularly varying function is WAP; in particular, if
	$f\in\mathrm{RV}_0(\beta)$ with $\beta>1$, Karamata's theorem yields
	$F\in\mathrm{RV}_0(1-\beta)$, $F^{-1}\in\mathrm{RV}_\infty(-1/(\beta-1))$, and
	$f\circ F^{-1}\in\mathrm{RV}_\infty(-\beta/(\beta-1))$, which is WAP. By
	contrast, when $\beta=1$, $F$ is slowly varying at $0$, so $F^{-1}$ and
	$f\circ F^{-1}$ are rapidly varying at $\infty$ and fail to be WAP; hence
	solutions decaying faster than every power law lie outside the WAP class.
	
	At the other extreme, rapidly varying (super-flat) functions $f$ at $0^+$ can
	still yield a WAP composition $\varphi = f\circ F^{-1}$. For example, when
	$f(x)=e^{-1/x}$, integration by parts gives
	\[
	F(x)\sim x^2 e^{1/x} = \frac{x^2}{f(x)} \qquad \text{as } x\to 0^+.
	\]
	Setting $x=F^{-1}(t)$ yields $(f\circ F^{-1})(t) \sim F^{-1}(t)^2/t$ as
	$t\to\infty$, while $\log F(x)\sim 1/x$ gives $F^{-1}(t)\sim 1/\log t$.
	Therefore,
	\[
	(f\circ F^{-1})(t)\sim \frac{1}{t(\log t)^2} \qquad \text{as } t\to\infty,
	\]
	so $f\circ F^{-1}\in\mathrm{RV}_\infty(-1)$ is WAP. Section~\ref{sec:sufficient_WAP}
	develops systematic derivative criteria showing that this behaviour holds
	across a wide class of weak nonlinearities without requiring explicit inversion
	of $F$.
	
	\subsection{Main characterisation theorems}
	
	We are now ready to state our main characterisations of rate-preserving
	perturbations.
	
	\begin{theorem} \label{thm.22}
		Let $f\in C([0,\infty);[0,\infty))$ be increasing with $f(0)=0$, $f(x)>0$ for
		$x>0$, and $f(x)/x\to0$ as $x\to0^+$. Suppose that $F$ is given by \eqref{eq:F}
		with $F(x)\to \infty$ as $x\to 0^+$, and that $f \circ F^{-1}$ is WAP at
		$\infty$. Let $g\in C([0,\infty);[0,\infty))$ and suppose $x$ is a $C^1$
		solution to
		\[
		x'(t) = -f(x(t)) + g(t), \qquad t \ge 0, \qquad x(0) = \xi > 0.
		\]
		Then the following are equivalent:
		\begin{itemize}
			\item[(A)] $\displaystyle \lim_{t \to \infty} \frac{F(x(t))}{t} = 1.$
			\item[(B)] $g \in L^1(0,\infty)$ and $\displaystyle \lim_{t \to \infty} \frac{1}{t} \int_0^t \frac{g(s)}{f\big(F^{-1}(s)\big)} \, ds = 0.$
			\item[(C)] $g \in L^1(0,\infty)$ and $\displaystyle \lim_{t \to \infty} \frac{1}{t (f \circ F^{-1})(t)} \int_t^{2t} g(s)\,ds = 0.$
		\end{itemize}
	\end{theorem}
	
	Under the WAP hypothesis on $f\circ F^{-1}$, the forward-window condition
	$\int_t^{2t} g(s)\,ds = o(t(f\circ F^{-1})(t))$ in~(C) is equivalent to the
	backward-window condition
	\[
	\lim_{t\to\infty} \frac{1}{t(f\circ F^{-1})(t)}\int_{t/2}^{t} g(s)\,ds = 0,
	\]
	which depends only on forcing values up to the present time $t$.
	
	For $f\in\mathrm{RV}_0(\beta)$ ($\beta>1$), Appleby and
	Patterson~\cite{appleby2014necessary} proved that rate preservation is
	characterised by the simpler tail-integral condition
	$\int_t^\infty g(s)\,ds = o(F^{-1}(t))$, and Appleby~\cite{appleby2017exact}
	extended this to asymptotic preserving functions. Our next theorem shows that
	this equivalence holds under the WAP hypothesis whenever $F$ satisfies a
	two-sided elasticity bound relative to $x/f(x)$.
	
	\begin{theorem} \label{thm.23}
		Let $f$, $F$, $g$, and $x$ satisfy the hypotheses of Theorem~\ref{thm.22}.
		Suppose further that
		\begin{equation} \label{eq.8}
			0 < \liminf_{x \to 0^+} \frac{F(x)}{x /f(x)}
			\le \limsup_{x \to 0^+} \frac{F(x)}{x/ f(x)} < +\infty.
		\end{equation}
		Then the following are equivalent:
		\begin{itemize}
			\item[(A)] $\displaystyle \lim_{t\to\infty} \frac{F(x(t))}{t} = 1.$
			\item[(B)] $g \in L^1(0,\infty)$ and
			\[
			\lim_{t\to\infty} \frac{1}{t} \int_0^t \frac{g(s)}{f(F^{-1}(s))} \, ds = 0.
			\]
			\item[(C)] $g \in L^1(0,\infty)$ and
			\[
			\lim_{t\to\infty} \frac{\int_t^\infty g(s)\,ds}{F^{-1}(t)} = 0.
			\]
		\end{itemize}
	\end{theorem}
	
	Structurally, the equivalence (A)~$\Longleftrightarrow$~(B) is obtained directly
	from the ODE (Section~\ref{sec:equiv_AB}), but condition~(B) couples $g$ and
	$f\circ F^{-1}$ inside the same integral and is inconvenient to verify when $g$
	is irregular. Decoupling~(B) into~(C) via the Ces\`aro--Tauberian theorems of
	Section~\ref{sec:cesaro_char} not only yields practical criteria, but also
	reflects the integral structure of \eqref{eq:main_ode}. Indeed, integrating
	\eqref{eq:main_ode} over $[0,t]$ gives
	$x(t) = \xi + \int_0^t g(s)\,ds - \int_0^t f(x(s))\,ds$: when
	$g\in L^1(0,\infty)$, the cumulative forcing shifts the effective initial
	condition by a finite constant. Setting $G(t) := \int_t^\infty g(s)\,ds$ and
	$\delta(t) := x(t)+G(t)$ yields
	\[
	\delta'(t) = -f(\delta(t)-G(t)), \qquad t\geq 0,
	\]
	showing that the trajectory $x = \delta - G$ depends on $g$ exclusively through
	its tail integral $G$. Conditions~(C) in Theorems~\ref{thm.22} and~\ref{thm.23}
	make this dependence explicit, since $\int_t^{2t} g(s)\,ds = G(t)-G(2t)$ and
	$\int_t^\infty g(s)\,ds = G(t)$.
	
	\subsection{Pointwise forcing and derivative asymptotics}
	
	Condition~(B) in Theorems~\ref{thm.22} and~\ref{thm.23} is clearly implied by
	the pointwise smallness condition
	\begin{equation} \label{eq.gisofFinv}
		\lim_{t\to\infty} \frac{g(t)}{(f\circ F^{-1})(t)}=0,
	\end{equation}
	because $f\circ F^{-1}\in L^1(0,\infty)$ ensures $g\in L^1(0,\infty)$ and a
	function tending to zero has zero Ces\`aro mean. While \eqref{eq.gisofFinv} is
	not necessary for $F(x(t))/t\to 1$, it precisely characterises simultaneous
	preservation of the decay rate of $x$ and the asymptotic behaviour of $x'$.
	
	\begin{theorem} \label{thm.24}
		Let $f$, $F$, $g$, and $x$ satisfy the hypotheses of Theorem~\ref{thm.22}.
		Then the following are equivalent:
		\begin{itemize}
			\item[(A)] $\displaystyle \lim_{t \to \infty} \frac{F(x(t))}{t} = 1$ and $\displaystyle x'(t)\sim -f(x(t))$ as $t\to\infty$.
			\item[(B)] $\displaystyle \lim_{t \to \infty} \frac{g(t)}{(f\circ F^{-1})(t)}=0.$
		\end{itemize}
	\end{theorem}
	
	\begin{proof}
		Suppose \eqref{eq.gisofFinv} holds. As noted above, condition~(B) of
		Theorem~\ref{thm.22} is satisfied, so $F(x(t))\sim t$ as $t\to\infty$. Since
		$f\circ F^{-1}$ is WAP at $\infty$, there exist $c>0$ and $T_1>0$ such that
		\[
		f(x(t))=(f\circ F^{-1})(F(x(t)))\geq c (f\circ F^{-1})(t), \qquad t\geq T_1.
		\]
		Combining this with \eqref{eq.gisofFinv} gives $g(t)=o(f(x(t)))$ as
		$t\to\infty$, whence $x'(t)=-f(x(t))+g(t)\sim -f(x(t))$. Thus~(B) implies~(A).
		
		Conversely, suppose~(A) holds. From $x'(t)\sim -f(x(t))$ and
		\eqref{eq:main_ode}, we have $g(t)=o(f(x(t)))$ as $t\to\infty$. Since
		$f\circ F^{-1}$ is WAP at $\infty$ and $F(x(t))\sim t$, there exist $C>0$ and
		$T_2>0$ such that
		\[
		f(x(t))=(f\circ F^{-1})(F(x(t)))\le C (f\circ F^{-1})(t), \qquad t\geq T_2.
		\]
		Hence
		\[
		0\leq \frac{g(t)}{(f\circ F^{-1})(t)}=\frac{g(t)}{f(x(t))}\cdot \frac{f(x(t))}{(f\circ F^{-1})(t)} \longrightarrow 0
		\qquad \text{as } t\to\infty,
		\]
		which is~(B). (Note also that $x'(t)\sim -f(x(t))$ alone already implies
		$F(x(t))\sim t$ by asymptotic integration.)
	\end{proof}
	
	Theorem~\ref{thm.24} can be sharpened when $\varphi = f\circ F^{-1}$ is
	\emph{strongly asymptotic preserving} (SAP) at $\infty$, meaning that $\varphi$
	is WAP at $\infty$ and the dilation bounds in \eqref{eq.WAPlimitsInf} satisfy
	$\lim_{\lambda\to 1} \underline{\Phi}_\lambda = \lim_{\lambda\to 1} \overline{\Phi}_\lambda = 1$.
	In that case, $a(t)\sim b(t)\to\infty$ implies $\varphi(a(t))\sim\varphi(b(t))$,
	so $F(x(t))\sim t$ yields $f(x(t))=(f\circ F^{-1})(F(x(t)))\sim (f\circ F^{-1})(t)$.
	
	\begin{theorem} \label{thm.24bis}
		Let $f$, $F$, $g$, and $x$ satisfy the hypotheses of Theorem~\ref{thm.22}, and
		suppose further that $f\circ F^{-1}$ is SAP at $\infty$. Then the following
		are equivalent:
		\begin{itemize}
			\item[(A)] $\displaystyle \lim_{t \to \infty} \frac{F(x(t))}{t} = 1$ and $\displaystyle x'(t)\sim -f(F^{-1}(t))$ as $t\to\infty$.
			\item[(B)] $\displaystyle \lim_{t \to \infty} \frac{g(t)}{(f\circ F^{-1})(t)}=0.$
		\end{itemize}
	\end{theorem}
	
	\begin{proof}
		If~(B) holds, Theorem~\ref{thm.24} gives $F(x(t))\sim t$ and
		$x'(t)\sim -f(x(t))$, and the SAP property then yields
		$-x'(t)\sim (f\circ F^{-1})(F(x(t)))\sim (f\circ F^{-1})(t)$. Conversely, if~(A)
		holds, the SAP property gives $-x'(t)\sim (f\circ F^{-1})(t)\sim (f\circ F^{-1})(F(x(t)))=f(x(t))$,
		so~(A) of Theorem~\ref{thm.24} holds and hence~(B) follows.
	\end{proof}
	
	\section{Optimality and the boundary of the characterisation}
	\label{sec:sharpness_counterexamples}
	
	Before proving Theorems~\ref{thm.22} and~\ref{thm.23}, we delineate the precise
	optimality boundary of their hypotheses. Our main theorems rest on three
	structural assumptions: sublinearity of $f$ at $0^+$, the WAP property for
	$\varphi = f\circ F^{-1}$ at $\infty$, and (in Theorem~\ref{thm.23}) the
	two-sided elasticity bound \eqref{eq.8}. We now summarise the explicit
	constructions proved in Appendices~\ref{app:rv1}--\ref{app:superflat}, which
	show that relaxing any one of these assumptions causes the characterising
	perturbation criteria to break down or change form fundamentally. Throughout
	this section we write $G(t):=\int_t^\infty g(s)\,ds$.
	
	\subsection{Sublinear $f\in\mathrm{RV}_0(1)$ with $f\circ F^{-1}$ non-WAP (Appendix~\ref{app:rv1})}
	
	When $f$ remains sublinear ($f(x)/x\to 0$) but lies in the borderline class
	$\mathrm{RV}_0(1)$ (satisfying the $C^2$ regularity conditions of
	Hypothesis~\ref{hyp:f_properties}, such as $f(x)=x/\log(1/x)$), $F$ is slowly
	varying at $0^+$ and both $F^{-1}$ and $f\circ F^{-1}$ are rapidly varying at
	$\infty$, hence non-WAP. In Theorems~\ref{thm:direct_alpha_construction}
	and~\ref{thm:direct_trajectory_properties} (Appendix~\ref{app:rv1}), we
	construct for every such $f$ an explicit positive, integrable forcing $g$
	(generated via a self-regulating shift ODE $\alpha'(t)=\exp(-\frac{1}{2}\alpha(t)f'(F^{-1}(t)))$)
	whose solution $x(t)=F^{-1}(t-\alpha(t))$ satisfies
	\[
	\lim_{t\to\infty}\frac{F(x(t))}{t}=1,
	\qquad \text{yet} \qquad
	\lim_{t\to\infty}\frac{g(t)}{(f\circ F^{-1})(t)}=\infty
	\quad \text{and} \quad
	\lim_{t\to\infty}\frac{G(t)}{F^{-1}(t)}=\infty.
	\]
	In particular, both the Ces\`aro condition~(B) of Theorem~\ref{thm.22} and the
	tail condition~(C) of Theorem~\ref{thm.23} fail to be necessary for rate
	preservation once $f\circ F^{-1}$ loses the WAP property.
	
	\subsection{Asymptotically linear $f$ with $f\circ F^{-1}$ non-WAP (Appendix~\ref{app:linear})}
	
	If sublinearity is relaxed to asymptotic linearity, $f(x)\sim ax$ as $x\to0^+$
	with $a>0$, then $F(x)\sim \frac{1}{a}\log(1/x)$ is slowly varying at $0^+$ and
	$f\circ F^{-1}$ is rapidly varying at $\infty$ (Remark~\ref{rem:linear_nonwap}).
	In this regime, Theorems~\ref{thm:linear_suff} and~\ref{thm:linear_nec}
	(Appendix~\ref{app:linear}) establish a \emph{different} exact necessary and
	sufficient characterisation of rate preservation:
	\[
	\lim_{t\to\infty}\frac{F(x(t))}{t}=1
	\quad \Longleftrightarrow \quad
	\lim_{t\to\infty}\frac{\log x(t)}{t}=-a
	\quad \Longleftrightarrow \quad
	g\in L^1(0,\infty) \text{ and } \liminf_{t\to\infty}\frac{F(G(t))}{t}\ge 1.
	\]
	Furthermore, Theorem~\ref{thm:linear_supercrit} (Appendix~\ref{app:linear})
	constructs an ultimately positive continuous forcing $g\in L^1(0,\infty)$ for
	which $\lim_{t\to\infty} F(G(t))/t=1$ and $F(x(t))/t\to 1$, while
	\[
	\lim_{t\to\infty}\frac{G(t)}{F^{-1}(t)}=\infty.
	\]
	Thus the tail criterion $G=o(F^{-1})$ is strictly stronger than necessary in
	the linearisable regime and is replaced by the logarithmic-scale condition
	$\liminf_{t\to\infty} F(G(t))/t \ge 1$.
	
	\subsection{Superlinear $f$ (Appendix~\ref{app:superlinear})}
	
	If $f$ is superlinear at zero in the sense that $f(x)/x\to\infty$ as $x\to0^+$
	while $F(x)\to\infty$, Theorem~\ref{thm:superlinear_nonwap}
	(Appendix~\ref{app:superlinear}) proves that $\varphi = f\circ F^{-1}$ is
	automatically rapidly varying at $\infty$, and hence never WAP. Moreover, when
	$f'(x)\to\infty$ as $x\to 0^+$, Theorem~\ref{thm:superlinear_counterex}
	(Appendix~\ref{app:superlinear}) constructs a positive continuous forcing
	$g\in L^1(0,\infty)$ (using the shift $\alpha(t)=(1+t)^{1/2}-1$) such that
	\[
	\lim_{t\to\infty}\frac{F(x(t))}{t}=1
	\quad \text{and} \quad
	\liminf_{t\to\infty}\frac{F(G(t))}{t}\ge 1,
	\qquad \text{yet} \qquad
	\lim_{t\to\infty}\frac{G(t)}{F^{-1}(t)}=\infty.
	\]
	Hence the sublinear WAP characterisation breaks down in the superlinear regime
	as well, and the same pattern $\liminf_{t\to\infty} F(G(t))/t\ge 1$ re-emerges.
	
	\subsection{Super-flat $f$ with $f\circ F^{-1}$ WAP: failure of the tail criterion (Appendix~\ref{app:superflat})}
	
	Finally, we test the sharpness of the two-sided elasticity hypothesis
	\eqref{eq.8} in Theorem~\ref{thm.23}. When $f$ is super-flat---for instance
	$f(x)=e^{-1/x}$---the composition $f\circ F^{-1}\in \mathrm{RV}_\infty(-1)$
	\emph{is} WAP, so Theorem~\ref{thm.22} applies and characterises rate
	preservation by the finite-window condition
	$\int_t^{2t} g(s)\,ds = o(t(f\circ F^{-1})(t))$. However, by
	Theorem~\ref{thm:rapid_variation_WAP}, $F(x) = o(x/f(x))$ as $x\to 0^+$, so
	\eqref{eq.8} fails.
	
	In Appendix~\ref{app:superflat} (Theorem~\ref{thm:stalling_counterex}), we
	prove that the weaker tail condition $G(t)=o(F^{-1}(t))$ is \emph{not
		sufficient} for rate preservation when $f(x)=e^{-1/x}$: for any prescribed
	$c\in(0,1)$, we construct a $C^1$ ``stalling'' trajectory
	$x(t)=F^{-1}(H(t))$ with nonnegative continuous forcing $g\in L^1(0,\infty)$
	such that
	\[
	\lim_{t\to\infty}\frac{G(t)}{F^{-1}(t)}=0,
	\qquad \text{yet} \qquad
	\liminf_{t\to\infty}\frac{F(x(t))}{t}=c < 1 = \limsup_{t\to\infty}\frac{F(x(t))}{t}.
	\]
	Thus, in the super-flat WAP regime, the full strength of the finite-window
	criterion in Theorem~\ref{thm.22}(C) cannot be weakened to the tail criterion
	of Theorem~\ref{thm.23}(C). For completeness, Appendix~\ref{app:superflat} also
	contains Proposition~\ref{prop:solution-asymptotic-not-tail}, showing that for
	the super-flat map $f(x)=x^2 e^{-1/x}$ (where $\varphi\in\mathrm{RV}_\infty(-1)$
	is WAP), there exists a positive forcing $g\notin L^1(0,\infty)$ with
	$x(t)\sim F^{-1}(t)$, underscoring why $F(x(t))/t\to 1$ is a cleaner
	rate-preservation criterion than $x(t)\sim F^{-1}(t)$.
	
	\subsection{Summary}
	
	Taken together, Appendices~\ref{app:rv1}--\ref{app:superflat} show that
	Theorems~\ref{thm.22} and~\ref{thm.23} are sharply formulated: removing WAP,
	sublinearity, or the two-sided elasticity bound \eqref{eq.8} invalidates the
	corresponding characterisations. Across the three non-WAP regimes
	(Appendices~\ref{app:rv1}--\ref{app:superlinear}), the evidence strongly
	suggests that the replacement characterisation when $f\circ F^{-1}$ is no
	longer WAP may take the form $\liminf_{t\to\infty} F(G(t))/t \ge 1$; a
	systematic treatment of that non-WAP regime is left for future work.
	
	\section{Ces\`aro and Tauberian characterisations}
	\label{sec:cesaro_char}
	
	In this section we establish two general real-variable equivalences linking
	Ces\`aro means of derivative ratios $a'/b'$ to dyadic increments of $a$ and to
	the quotient $a/b$. These results underpin the decoupling of condition~(B) into
	condition~(C) in Theorems~\ref{thm.22} and~\ref{thm.23} (via $a(t)=\int_t^\infty g(s)\,ds$
	and $b(t)=F^{-1}(t)$), and are also of independent interest as Tauberian-type
	converse principles to l'H\^opital's rule.
	
	\begin{theorem} \label{thm.generalcesaro}
		Suppose $a \in C^1([0,\infty);[0,\infty))$ is non-increasing, and
		$b \in C^1([0,\infty);(0,\infty))$ is decreasing with $b'(t)<0$ for all
		$t\ge 0$ and $-b'$ WAP at $\infty$.
		Then the following are equivalent:
		\begin{itemize}
			\item[(A)]
			\[
			\lim_{t \to \infty} \frac{1}{t} \int_0^t \frac{a'(s)}{b'(s)}\,ds = 0,
			\]
			\item[(B)]
			\[
			\lim_{t \to \infty} \frac{a(t) - a(2t)}{t |b'(t)|} = 0.
			\]
		\end{itemize}
	\end{theorem}
	
	\begin{proof}
		We start by proving that $(B) \Longrightarrow (A)$. Define
		\[
		C(t) := \frac{1}{t} \int_0^t \left| \frac{a'(s)}{b'(s)} \right| ds
		= \frac{1}{t} \int_0^t \frac{a'(s)}{b'(s)}\,ds.
		\]
		Note that if $C(2^N) \to 0$ as $N \to \infty$ through integers, then $C(t) \to 0$
		as $t \to \infty$ through the reals (since for $2^N \le t < 2^{N+1}$,
		$0 \le C(t) \le 2\,C(2^{N+1})$). Define
		\[
		\eta(t) := \frac{a(t) - a(2t)}{t |b'(t)|} \geq 0,
		\]
		and let $\varphi(t):= -b'(t) = |b'(t)|$, so $\varphi$ is WAP at $\infty$ and
		eventually non-increasing. For $s \in [2^K, 2^{K+1}]$ (with $K$ sufficiently
		large), monotonicity gives $\varphi(2^K) \ge \varphi(s) \ge \varphi(2^{K+1})$, so
		\[
		\frac{1}{\varphi(s)} \le \frac{1}{\varphi(2^{K+1})}
		= \frac{1}{\varphi(2^K)} \cdot \frac{\varphi(2^K)}{\varphi(2^{K+1})}.
		\]
		Since $\varphi$ is WAP at $\infty$, there exists a finite $M>0$ such that
		\[
		\limsup_{t \to \infty} \frac{\varphi(t/2)}{\varphi(t)} =: M < +\infty.
		\]
		Therefore, there exists $T_0 > 0$ such that $\varphi(t/2)/\varphi(t) \le 2M$
		for all $t \ge T_0$ (and $\varphi$ is non-increasing on $[T_0/2,\infty)$).
		Choose $K_0 \in \mathbb{N}$ such that $2^{K_0 + 1} \ge T_0$. Then for every
		integer $K \ge K_0$ and $s \in [2^K, 2^{K+1}]$, we have
		\begin{equation} \label{eq.upper1overphi}
			\frac{1}{\varphi(s)} \le \frac{1}{\varphi(2^K)} \cdot 2M.
		\end{equation}
		Thus, for $K \ge K_0$, using \eqref{eq.upper1overphi}, the monotonicity of $a$,
		and the definition of $\eta$, we obtain
		\begin{align} \label{eq.a'b'dyadicint}
			\int_{2^K}^{2^{K+1}} \frac{|a'(s)|}{|b'(s)|} \, ds
			&= \int_{2^K}^{2^{K+1}} \frac{|a'(s)|}{\varphi(s)} \, ds \nonumber \\
			&\le 2M \cdot \frac{1}{\varphi(2^K)} \int_{2^K}^{2^{K+1}} |a'(s)| \, ds \nonumber \\
			&= 2M \cdot 2^K \cdot \eta(2^K) = 2M \cdot 2^K \cdot \tilde{\eta}_K,
		\end{align}
		where $\tilde{\eta}_K := \eta(2^K) \to 0$ as $K\to \infty$ by~(B). Now let
		$N \ge K_0 + 1$. Using the definition of $C(t)$ and \eqref{eq.a'b'dyadicint},
		we get
		\[
		C(2^N) \le \frac{2^{K_0}}{2^N} C(2^{K_0})
		+ \frac{1}{2^N} \sum_{K=K_0}^{N-1} 2M \cdot 2^K \tilde{\eta}_K
		= \frac{2^{K_0}}{2^N} C(2^{K_0}) + 2M \cdot \frac{1}{2^N} \sum_{K=K_0}^{N-1} 2^K \tilde{\eta}_K.
		\]
		Since $\tilde{\eta}_K \to 0$ as $K \to \infty$, the Stolz--Ces\`aro lemma
		ensures that the last term on the right-hand side tends to zero as
		$N\to\infty$. Likewise, since $2^{K_0}C(2^{K_0})$ is finite, the first term
		tends to zero as $N\to\infty$. Hence $C(2^N)\to 0$ as $N\to\infty$, completing
		the proof that (B)~implies~(A).
		
		Now we show that $(A) \Longrightarrow (B)$. Define
		\[
		A(t) = \int_1^t \frac{|a'(s)|}{|b'(s)|} \, ds.
		\]
		Since $C(t) \to 0$ as $t \to \infty$, we have $A(t)=o(t)$ as $t\to\infty$.
		Because $\varphi=|b'|$ is eventually non-increasing, for sufficiently large $t$
		and $s \in [t,2t]$ we have $\varphi(2t) \le \varphi(s) \le \varphi(t)$, and
		thus
		\[
		\frac{1}{|b'(s)|}=\frac{1}{\varphi(s)} \ge \frac{1}{\varphi(t)} = \frac{1}{|b'(t)|}.
		\]
		Using this inequality, the definition of $A$, and the monotonicity of $a$, we
		obtain
		\[
		2 \cdot \frac{A(2t)}{2t} - \frac{A(t)}{t} = \frac{A(2t) - A(t)}{t}
		= \frac{1}{t} \int_t^{2t} \frac{|a'(s)|}{|b'(s)|}\, ds
		\ge \frac{a(t) - a(2t)}{t|b'(t)|} \ge 0.
		\]
		Since $A(t)=o(t)$ as $t\to\infty$, the left-hand side tends to zero as
		$t\to\infty$, and therefore so does the right-hand side. This proves~(B).
	\end{proof}
	
	Condition~(B) in Theorem~\ref{thm.generalcesaro} decouples the pointwise
	integrand in~(A). With a two-sided elasticity bound on $b$, both conditions
	become equivalent to the direct quotient limit $a(t)/b(t)\to 0$. Since the
	direct converse to l'H\^opital's rule ($a/b\to 0 \implies a'/b'\to 0$) fails in
	general, the following theorem provides a sharp Ces\`aro-averaged converse.
	
	\begin{theorem} \label{thm.converselhopital}
		Suppose $a\in C^1([0,\infty);[0,\infty))$ is non-increasing with $a(t) \to 0$
		as $t\to\infty$. Suppose further that $b\in C^1([0,\infty);(0,\infty))$ is
		decreasing with $b'(t)<0$ for all $t\ge 0$ and $-b'$ WAP at $\infty$, and
		assume there exist constants $c, C > 0$ such that
		\begin{equation} \label{eq.bvonmises}
			c \le \frac{t|b'(t)|}{b(t)} \le C, \qquad t \ge 1.
		\end{equation}
		Then the condition
		\[
		\textnormal{(C)} \qquad \lim_{t \to \infty} \frac{a(t)}{b(t)} = 0
		\]
		is equivalent to conditions~\textnormal{(A)} and~\textnormal{(B)} of
		Theorem~\ref{thm.generalcesaro}.
	\end{theorem}
	
	\begin{proof}
		We first show that (C)~implies~(B). From \eqref{eq.bvonmises} and~(C),
		\[
		\frac{a(t)}{t\, |b'(t)|} = \frac{a(t)}{b(t)}\cdot
		\frac{b(t)}{t|b'(t)|}\to 0 \qquad \text{as } t\to\infty.
		\]
		Since $b$ is decreasing, $b(2t)/b(t)<1$, so likewise
		\[
		\frac{a(2t)}{t\, |b'(t)|} = \frac{a(2t)}{b(2t)}\cdot\frac{b(2t)}{b(t)} \cdot
		\frac{b(t)}{t|b'(t)|}\to 0 \qquad \text{as } t\to\infty.
		\]
		Taking the difference of these two limits yields~(B).
		
		Next, we show that (A)~implies~(C), which completes the equivalence since
		(C)~$\Rightarrow$~(B) and (B)~$\Leftrightarrow$~(A) by
		Theorem~\ref{thm.generalcesaro}. First, for $t \ge 1$ and any integer $K \ge 0$,
		condition \eqref{eq.bvonmises} yields
		\[
		- \log b(2^K t) + \log b(t)
		= \int_t^{2^K t} \frac{-b'(s)}{b(s)} \, ds
		\ge \int_t^{2^K t} \frac{c}{s} \, ds
		= c \log(2^K),
		\]
		which rearranges to
		\begin{equation} \label{eq.bdyadic}
			\frac{b(2^K t)}{b(t)} \le (2^{-c})^K,
			\qquad t \ge 1, \quad K \ge 0.
		\end{equation}
		Fix $t \ge 1$ and define $j_t(K):= a(2^K t) - a(2^{K+1} t) \ge 0$ for integers
		$K\ge 0$. For every $N \ge 1$, $\sum_{K=0}^{N-1} j_t(K) = a(t) - a(2^N t)$.
		Since $a(2^N t)\to 0$ as $N\to\infty$, letting $N\to\infty$ gives
		\begin{equation} \label{eq.sumrep}
			a(t)=\sum_{K=0}^{\infty} j_t(K) = \sum_{K=0}^{\infty} \left[ a(2^K t) - a(2^{K+1} t) \right].
		\end{equation}
		Assuming~(A), Theorem~\ref{thm.generalcesaro} gives~(B), so for every
		$\varepsilon > 0$ there exists $T(\varepsilon) > 1$ such that
		$a(u) - a(2u) \le \varepsilon\, u |b'(u)| \le \varepsilon C b(u)$ for all
		$u\ge T(\varepsilon)$, using \eqref{eq.bvonmises} in the last inequality. For
		$t \ge T(\varepsilon)$ and every $K\ge 0$, we have $2^K t \ge t \ge T(\varepsilon)$,
		so $a(2^K t) - a(2^{K+1} t) \le \varepsilon C\, b(2^K t)$. Applying
		\eqref{eq.sumrep}, this bound, and \eqref{eq.bdyadic}, we obtain
		\[
		0\le \frac{a(t)}{b(t)} \le \varepsilon C \sum_{K=0}^{\infty} \frac{b(2^K t)}{b(t)}
		\le \varepsilon C \sum_{K=0}^\infty (2^{-c})^K =: \varepsilon C',
		\qquad t \ge T(\varepsilon)>1,
		\]
		where $C'<\infty$ is independent of $\varepsilon$. Since $\varepsilon>0$ is
		arbitrary, $a(t)/b(t)\to 0$ as $t\to\infty$, proving~(C).
	\end{proof}
	
	\section{Proof of the main ODE characterisations}
	\label{sec:equiv_AB}
	
	\subsection{The coupled Ces\`aro characterisation}
	
	Statements~(A) and~(B) in Theorems~\ref{thm.22} and~\ref{thm.23} are identical.
	In this section, we establish their equivalence directly from the differential
	equation \eqref{eq:main_ode} under the hypotheses of Theorem~\ref{thm.22}.
	
	\begin{proof}[Proof of \textnormal{(A)} $\Longleftrightarrow$ \textnormal{(B)} in Theorems~\ref{thm.22} and~\ref{thm.23}]
		We start by proving $(B) \Longrightarrow (A)$. Define $y$ as the solution of
		\[
		y'(t) = -f(y(t)), \qquad t \ge 0, \qquad y(0) = \xi' = \frac{1}{2}\min(\xi,1).
		\]
		Since $g(t)\ge 0$ and $x(0)=\xi > \xi'$, comparison gives $x(t) \ge y(t)$ for
		all $t\geq 0$, and hence $F(x(t)) \le F(y(t))$ for all $t\geq 0$. Because
		$F(y(t)) - F(\xi') = t$, we obtain
		\begin{equation} \label{eq.Fxupper}
			F(x(t)) \le F(\xi') + t, \qquad t \ge 0.
		\end{equation}
		Dividing \eqref{eq:main_ode} by $f(x(t))$ gives
		\[
		-\frac{x'(t)}{f(x(t))} = 1 - \frac{g(t)}{f(x(t))}.
		\]
		Integrating over $[0,t]$ and using the definition of $F$ yields
		$F(x(t)) - F(\xi) = t - \int_0^t \frac{g(s)}{f(x(s))}\,ds$, or
		\begin{equation} \label{eq.Fxtmaster}
			\frac{F(x(t))}{t} = 1 + \frac{F(\xi)}{t} - \frac{1}{t} \int_0^t \frac{g(s)}{f(x(s))}\,ds.
		\end{equation}
		It therefore suffices to show that the last term on the right-hand side tends
		to zero. From \eqref{eq.Fxupper}, $x(t) \ge F^{-1}(F(\xi')+t)$ for all
		$t\geq 0$. Since $f$ is increasing and $g\ge 0$, integrating over $[0,t]$ gives
		\begin{equation} \label{eq.gaveest1}
			0\le \frac{1}{t} \int_0^t \frac{g(s)}{f(x(s))}\,ds
			\le \frac{1}{t} \int_0^t \frac{g(s)}{(f \circ F^{-1})(F(\xi')+s)}\,ds.
		\end{equation}
		Because $\xi' < 1$, we have $F(\xi') > 0$. For $s \ge F(\xi')$, we have
		$2s \ge s + F(\xi')$, and since $f \circ F^{-1}$ is decreasing,
		$(f \circ F^{-1})(2s) \le (f \circ F^{-1})(s + F(\xi'))$. Because
		$f \circ F^{-1}$ is WAP at $\infty$, there exists $C>0$ such that
		\[
		\liminf_{t \to \infty} \frac{(f \circ F^{-1})(2t)}{(f \circ F^{-1})(t)} = C > 0,
		\]
		so there exists $T_1>0$ such that $(f \circ F^{-1})(2t)/(f \circ F^{-1})(t) \ge C/2$
		for all $t \ge T_1$. Setting $T_2 := \max(T_1, F(\xi'))$, we obtain for all
		$s \ge T_2$,
		\[
		\frac{g(s)}{(f \circ F^{-1})(s+F(\xi'))}
		\le \frac{g(s)}{(f \circ F^{-1})(2s)}
		= \frac{g(s)}{(f \circ F^{-1})(s)} \cdot \frac{(f \circ F^{-1})(s)}{(f \circ F^{-1})(2s)}
		\le \frac{2}{C} \cdot \frac{g(s)}{(f \circ F^{-1})(s)}.
		\]
		Splitting the integral on the right-hand side of \eqref{eq.gaveest1} at $T_2$
		gives, for all $t \ge T_2$,
		\[
		0\le \frac{1}{t} \int_0^t \frac{g(s)}{f(x(s))}\,ds
		\le \frac{1}{t} \int_0^{T_2} \frac{g(s)}{(f \circ F^{-1})(s+F(\xi'))}\,ds
		+ \frac{2}{C}\frac{1}{t} \int_{T_2}^t \frac{g(s)}{(f \circ F^{-1})(s)}\,ds.
		\]
		The first term tends to zero trivially as $t\to\infty$, and the second tends to
		zero by~(B) and the nonnegativity of $g$. Returning to \eqref{eq.Fxtmaster}
		yields $F(x(t))/t\to 1$ as $t\to\infty$, proving $(B) \Longrightarrow (A)$.
		
		Now we prove $(A) \Longrightarrow (B)$. From \eqref{eq.Fxtmaster} and~(A), we
		immediately have
		\begin{equation} \label{eq.gfint1}
			\lim_{t \to \infty} \frac{1}{t} \int_0^t \frac{g(s)}{f(x(s))} \, ds = 0.
		\end{equation}
		Moreover, (A)~implies that for any $\varepsilon \in (0,1)$, there exists
		$T(\varepsilon)>0$ such that $(1 - \varepsilon)t < F(x(t))$ for all
		$t \geq T(\varepsilon)$. Since $f \circ F^{-1}$ is decreasing and $g\ge 0$,
		\[
		0 \le \frac{g(t)}{f(F^{-1}((1-\varepsilon)t))}
		\le \frac{g(t)}{f(x(t))}, \qquad t \geq T(\varepsilon),
		\]
		and hence \eqref{eq.gfint1} yields
		\begin{equation} \label{eq.gfint2}
			\lim_{t \to \infty} \frac{1}{t} \int_0^t \frac{g(s)}{(f\circ F^{-1})((1-\varepsilon)s)} \, ds = 0.
		\end{equation}
		Let $\varphi = f \circ F^{-1}$. Since $\varphi$ is WAP at $\infty$, for each
		fixed $\varepsilon\in(0,1)$ there exist $C_\varepsilon<\infty$ and
		$s(\varepsilon)>0$ such that $\varphi((1-\varepsilon)s)\le C_\varepsilon\varphi(s)$
		for all $s>s(\varepsilon)$. Consequently, for $s>s(\varepsilon)$,
		\begin{equation} \label{eq.fFinvwap}
			\frac{g(s)}{\varphi(s)}
			\le C_\varepsilon \frac{g(s)}{\varphi((1-\varepsilon)s)}
			= C_\varepsilon \frac{g(s)}{f(F^{-1}((1-\varepsilon)s))}.
		\end{equation}
		Combining \eqref{eq.gfint2} and \eqref{eq.fFinvwap} yields
		\[
		\lim_{t\to\infty} \frac{1}{t} \int_0^t \frac{g(s)}{f(F^{-1}(s))}\,ds = 0,
		\]
		which is the second part of~(B).
		
		It remains to show that (A)~implies $g \in L^1(0,\infty)$. Integrating
		\eqref{eq:main_ode} over $[0,t]$ gives
		\begin{equation} \label{eq.gint}
			\int_0^t g(s)\,ds = \int_0^t f(x(s))\,ds + x(t) - \xi, \qquad t\geq 0.
		\end{equation}
		Since $f \circ F^{-1}$ is decreasing and $F(x(t))>(1-\varepsilon)t$ for
		$t\ge T(\varepsilon)$, we have $f(x(t)) \le (f \circ F^{-1})((1-\varepsilon)t)$
		for $t \geq T(\varepsilon)$. Substituting $u=F^{-1}((1-\varepsilon)s)$ gives
		\[
		\int_{T(\varepsilon)}^t f(x(s))\,ds
		\le \int_{T(\varepsilon)}^t (f \circ F^{-1})((1-\varepsilon)s)\,ds
		=\frac{1}{1-\varepsilon}\Big(F^{-1}((1-\varepsilon)T(\varepsilon))-F^{-1}((1-\varepsilon)t)\Big),
		\]
		which is uniformly bounded above as $t\to\infty$. Thus $t \mapsto f(x(t))$
		belongs to $L^1(0,\infty)$. Since $g\ge 0$ and $x(t)\to 0$ as $t\to\infty$,
		\eqref{eq.gint} shows that $\int_0^t g(s)\,ds$ converges to a finite limit, so
		$g\in L^1(0,\infty)$. Hence (A)~implies~(B), completing the proof.
	\end{proof}
	
	\subsection{Decoupling the forcing conditions}
	\label{sec:apply_cesaro}
	
	We now combine the Ces\`aro--Tauberian equivalences of
	Section~\ref{sec:cesaro_char} with the coupled characterisation of
	Section~\ref{sec:equiv_AB} to complete the proofs of Theorems~\ref{thm.22}
	and~\ref{thm.23}.
	
	Set $b(t):=F^{-1}(t)$. Then $b\in C^1([0,\infty);(0,\infty))$ is decreasing and
	$-b'(t)=(f\circ F^{-1})(t)>0$ for all $t\ge 0$ is WAP at $\infty$, so $b$
	satisfies all hypotheses of Theorem~\ref{thm.generalcesaro}. Under the
	additional elasticity condition \eqref{eq.8} of Theorem~\ref{thm.23}, the
	identity
	\[
	\frac{F(x)}{x /f(x)} \Big|_{x = F^{-1}(t)}
	= \frac{t}{F^{-1}(t)/ f(F^{-1}(t))}
	= \frac{t |b'(t)|}{b(t)}
	\]
	shows that $b$ also satisfies \eqref{eq.bvonmises} (on a half-line $[T,\infty)$,
	and hence on $[1,\infty)$ by continuity and positivity), so all hypotheses of
	Theorem~\ref{thm.converselhopital} hold.
	
	\begin{proof}[Proof of Theorem~\ref{thm.22}]
		By Section~\ref{sec:equiv_AB}, (A)~and~(B) are equivalent; it suffices to prove
		that (B)~and~(C) are equivalent. Both~(B) and~(C) include $g\in L^1(0,\infty)$,
		under which we define $a(t):=\int_t^\infty g(s)\,ds$. Then
		$a\in C^1([0,\infty);[0,\infty))$ is non-increasing with $-a'(t)=g(t)$ and
		$a(t)-a(2t)=\int_t^{2t} g(s)\,ds$. With this choice of $a$ and $b$, the
		Ces\`aro limit in Theorem~\ref{thm.22}(B) is precisely
		Theorem~\ref{thm.generalcesaro}(A), while the finite-window limit in
		Theorem~\ref{thm.22}(C) is precisely Theorem~\ref{thm.generalcesaro}(B).
		Applying Theorem~\ref{thm.generalcesaro} therefore establishes
		(B)~$\Longleftrightarrow$~(C) and completes the proof.
	\end{proof}
	
	\begin{proof}[Proof of Theorem~\ref{thm.23}]
		Again, (A)~and~(B) are equivalent by Section~\ref{sec:equiv_AB}, so it suffices
		to prove (B)~$\Longleftrightarrow$~(C). Under $g\in L^1(0,\infty)$, define
		$a(t):=\int_t^\infty g(s)\,ds$; then $a(t)\to 0$ as $t\to\infty$, and both $a$
		and $b=F^{-1}$ satisfy all hypotheses of Theorem~\ref{thm.converselhopital}.
		The Ces\`aro limit in Theorem~\ref{thm.23}(B) coincides with
		Theorem~\ref{thm.generalcesaro}(A), while the tail-integral limit in
		Theorem~\ref{thm.23}(C) is $a(t)/b(t)\to 0$, which is
		Theorem~\ref{thm.converselhopital}(C). Applying Theorem~\ref{thm.converselhopital}
		gives (B)~$\Longleftrightarrow$~(C) and completes the proof.
	\end{proof}
	
	\section{Checkable conditions for the WAP property}
	\label{sec:sufficient_WAP}
	
	In applications, finding $F^{-1}$ and $f\circ F^{-1}$ explicitly can be
	delicate, since it requires quadrature, inversion, and composition. It is
	therefore highly desirable to have checkable sufficient conditions for
	$f\circ F^{-1}$ to be WAP expressed directly in terms of $f$ and its low-order
	derivatives.
	
	When $f\in C^2$ with $f''>0$ near zero---a natural assumption given $f'(x)>0$
	and $f'(0^+)=0$---the key diagnostic quantity is
	\begin{equation} \label{eq.Q}
		Q(x)=\frac{f(x)f''(x)}{(f'(x))^2}.
	\end{equation}
	We show below that the two-sided bound
	\begin{equation} \label{eq.Qcondn}
		0 < c = \liminf_{x\to 0^+} Q(x) \le \limsup_{x\to 0^+} Q(x) = C < +\infty
	\end{equation}
	is sufficient to ensure that $f\circ F^{-1}$ is WAP at $\infty$, and that it
	holds across a broad collection of weak nonlinearities yielding power-law or
	slower-than-power-law decay.
	
	\subsection{Nontrivial $Q$ bounds imply WAP}
	
	To build intuition, note that if $f$ is smoothly regularly varying at $0^+$ with
	index $\beta>0$, so that $xf'(x)/f(x)\to\beta$ and $xf''(x)/f'(x)\to\beta-1$,
	then
	\[
	Q(x)=\frac{f(x)}{xf'(x)}\cdot\frac{xf''(x)}{f'(x)}\longrightarrow \frac{\beta-1}{\beta} \qquad \text{as } x\to 0^+.
	\]
	Thus \eqref{eq.Qcondn} holds for every $\beta>1$, fails at $\beta=1$ (where
	$f\circ F^{-1}$ is rapidly varying and non-WAP), and persists in the limit
	$\beta\to\infty$ corresponding to rapidly varying super-flat nonlinearities.
	
	\begin{theorem} \label{thm:Q_implies_WAP}
		Suppose that $f\in C^2((0,\delta))$ for some $\delta>0$ is increasing with
		increasing derivative, $f(0^+)=0$, and $F(x)\to\infty$ as $x\to 0^+$. Suppose
		that $Q$ defined by \eqref{eq.Q} obeys \eqref{eq.Qcondn}. Then
		$\varphi:=f\circ F^{-1}$ is WAP at $\infty$. In fact, for every $\lambda>1$,
		\[
		\lambda^{-1/c}\le \liminf_{t\to\infty} \frac{\varphi(\lambda t)}{\varphi(t)}
		\le \limsup_{t\to\infty}\frac{\varphi(\lambda t)}{\varphi(t)} \le \lambda^{-1/C}.
		\]
	\end{theorem}
	
	\begin{proof}
		Write $b(t) = F^{-1}(t)$, so $\varphi(t) = f(b(t))$ and
		$\varphi'(t) = -f'(b(t))f(b(t))$ for sufficiently large $t$. Hence
		\[
		\frac{-t \varphi'(t)}{\varphi(t)} = t f'(b(t)) = \frac{F(x)}{1/f'(x)}\Bigg|_{x=b(t)}.
		\]
		Since $b(t)\to 0^+$ as $t\to\infty$, the asymptotic behaviour of
		$-t\varphi'(t)/\varphi(t)$ as $t\to\infty$ is governed by $F(x)/(1/f'(x))$ as
		$x\to 0^+$.
		
		By \eqref{eq.Qcondn}, for every $\epsilon\in (0,c)$ there exists
		$\eta(\epsilon)\in(0,\delta)$ such that $c-\epsilon < Q(u) < C+\epsilon$ for all
		$u\in(0,\eta)$. Multiplying by $1/f(u)$ and integrating over $[x,\eta]$ gives
		\[
		(c-\epsilon)\int_x^\eta \frac{1}{f(u)}\,du < \int_x^\eta \frac{f''(u)}{(f'(u))^2}\,du < (C+\epsilon)\int_x^\eta\frac{1}{f(u)}\,du, \qquad 0<x<\eta.
		\]
		Using the fundamental theorem of calculus, we write
		\[
		\frac{F(x)}{1/f'(x)}
		= \frac{F(\eta)+\int_x^\eta \frac{1}{f(u)}\,du}{\frac{1}{f'(\eta)}+\int_x^\eta \frac{f''(u)}{(f'(u))^2} \,du}
		< \frac{F(\eta)+\int_x^\eta \frac{1}{f(u)}\,du}{\frac{1}{f'(\eta)}+(c-\epsilon)\int_x^\eta \frac{1}{f(u)}\,du}.
		\]
		Since $F(x)\to\infty$ as $x\to 0^+$, letting $x\to 0^+$ and then
		$\epsilon\to 0^+$ yields $\limsup_{x\to 0^+} \frac{F(x)}{1/f'(x)}\le \frac{1}{c}$.
		The analogous lower bound gives $\liminf_{x\to 0^+} \frac{F(x)}{1/f'(x)}\ge \frac{1}{C}$.
		Therefore,
		\[
		\frac{1}{C}\le \liminf_{t\to\infty} \frac{-t\varphi'(t)}{\varphi(t)}\le
		\limsup_{t\to\infty} \frac{-t\varphi'(t)}{\varphi(t)}\le \frac{1}{c}.
		\]
		Hence for every $\epsilon>0$ there exists $T(\epsilon)>0$ such that
		\[
		\left(\frac{1}{C}-\epsilon\right)\frac{1}{s}<\frac{-\varphi'(s)}{\varphi(s)}<\left(\frac{1}{c}+\epsilon\right)\frac{1}{s}, \qquad s\ge T(\epsilon).
		\]
		For any $\lambda>1$ and $t\ge T(\epsilon)$, integrating over $[t,\lambda t]$
		yields
		\[
		\left(\frac{1}{C}-\epsilon\right)\log\lambda < \log\left(\frac{\varphi(t)}{\varphi(\lambda t)}\right) < \left(\frac{1}{c}+\epsilon\right)\log\lambda,
		\]
		which rearranges to
		\[
		\lambda^{-1/c-\epsilon}<\frac{\varphi(\lambda t)}{\varphi(t)}<\lambda^{-1/C+\epsilon}, \qquad t\ge T(\epsilon).
		\]
		Letting $t\to\infty$ and then $\epsilon\to 0^+$ yields the stated bounds for
		$\lambda>1$. The case $\lambda\in(0,1)$ follows identically (or by replacing
		$\lambda$ with $1/\lambda$ and inverting), so $\varphi$ is WAP at $\infty$.
	\end{proof}
	
	\subsection{Rapidly varying $f$ and WAP $\varphi$}
	
	We now make precise the class of rapidly varying (super-flat) functions
	satisfying $Q(x)\to 1$, and show that this limit forces the two-sided
	elasticity condition \eqref{eq.8} of Theorem~\ref{thm.23} to fail.
	
	\begin{theorem} \label{thm:exp_rho_Q1}
		Suppose that $\rho$ is a $C^2$ regularly varying function at $\infty$ with
		index $\alpha>0$ satisfying
		\[
		\lim_{t\to\infty} \frac{t\rho'(t)}{\rho(t)}=\alpha, \qquad
		\lim_{t\to\infty} \frac{t\rho''(t)}{\rho'(t)}=\alpha-1.
		\]
		Define $f(x)=\exp(-\rho(1/x))$. Then for sufficiently small $x>0$, $f(x)>0$,
		$f'(x)>0$, $f''(x)>0$, and $f(x)=o(x^n)$ as $x\to 0^+$ for every
		$n\in\mathbb{N}$. Moreover, $Q(x)\to 1$ as $x\to 0^+$, so
		$\varphi = f\circ F^{-1}$ is WAP at $\infty$.
	\end{theorem}
	
	Choosing $\rho(t)=t$ ($\alpha=1$) recovers the prototypical super-flat function
	$f(x)=e^{-1/x}$.
	
	\begin{proof}
		Since $\log f(x)=-\rho(1/x)$, we have
		\[
		\frac{f'(x)}{f(x)}=\rho'(1/x)\cdot \frac{1}{x^2}=:a_1(x),
		\]
		and differentiating again gives
		\[
		\frac{f(x)f''(x)-(f'(x))^2}{(f(x))^2}
		=-\rho''(1/x)\frac{1}{x^4}-\rho'(1/x)\frac{2}{x^3}.
		\]
		Hence
		\[
		\frac{f''(x)}{f(x)}=a_1^2(x)-\rho''(1/x)\frac{1}{x^4}-\rho'(1/x)\frac{2}{x^3}=:a_2(x),
		\]
		and therefore
		\[
		Q(x)=\frac{a_2(x)}{a_1^2(x)}
		=1-\frac{\rho''(1/x)}{(\rho'(1/x))^2}-2\frac{\rho(1/x)}{\rho'(1/x)(1/x)}\cdot \frac{1}{\rho(1/x)}.
		\]
		Since $\alpha>0$, we have $\rho(1/x)\to \infty$ as $x\to 0^+$, so the third
		term tends to $0$. For the second term, as $t\to\infty$,
		\[
		\frac{\rho''(t)}{(\rho'(t))^2}
		=\frac{t\rho''(t)}{\rho'(t)}\cdot \frac{\rho(t)}{t\rho'(t)}\cdot \frac{1}{\rho(t)}
		\longrightarrow (\alpha-1)\cdot \frac{1}{\alpha}\cdot 0 = 0.
		\]
		Hence $Q(x)\to 1$ as $x\to 0^+$. In particular $f''(x)>0$ for small $x>0$,
		Theorem~\ref{thm:Q_implies_WAP} ensures that $\varphi$ is WAP at $\infty$, and
		Lemma~\ref{lem:elasticity} below confirms $f(x)=o(x^n)$ for all
		$n\in\mathbb{N}$.
	\end{proof}
	
	\begin{lemma} \label{lem:elasticity}
		Let $f \in C^2((0, x_0))$ satisfy $f(x) > 0$ and $f'(x) > 0$ on $(0, x_0)$,
		with $\lim_{x \to 0^+} f(x) = 0$. If $Q(x)\to 1$ as $x\to 0^+$, then
		\[
		\lim_{x \to 0^+} \frac{x f'(x)}{f(x)} = \infty,
		\]
		and consequently $f(x)=o(x^n)$ as $x\to 0^+$ for every $n\in\mathbb{N}$.
	\end{lemma}
	
	\begin{proof}
		Set $k(x) := f(x)/f'(x) > 0$. Then $k'(x) = 1 - Q(x) \to 0$ as $x\to 0^+$.
		Hence $k$ is uniformly continuous on $(0, x_0/2]$ and approaches a finite limit
		$L\ge 0$ as $x\to 0^+$. If $L>0$, then $f'(x)/f(x)\to 1/L\in(0,\infty)$, so
		integrating from $x$ to $c\in(0,x_0)$ would force $\log f(x)$ to remain
		bounded as $x\to 0^+$, contradicting $f(0^+)=0$. Thus $k(0^+)=0$.
		
		Extending $k$ continuously by $k(0)=0$, for any $\epsilon>0$ there exists
		$\delta>0$ such that $|k'(s)|\le\epsilon$ for $0<s<\delta$. Integrating $k'$
		over $[r,x]$ and letting $r\downarrow 0$ gives $0<k(x)\le\epsilon x$ for
		$0<x<\delta$. Hence $k(x)/x\to 0$, and therefore $xf'(x)/f(x)=x/k(x)\to\infty$.
		
		Consequently, for any $M\ge 1$ there exists $X(M)\in(0,1)$ such that
		$f'(u)/f(u)>M/u$ for $u\in(0,X(M))$. Integrating over $[x,X(M)]$ gives
		\[
		\log\frac{f(X(M))}{f(x)}
		> M\log\frac{X(M)}{x}.
		\]
		Thus $\log(1/f(x))/\log(1/x)\to\infty$, which implies
		$f(x)=o(x^n)$ for every $n\in\mathbb{N}$.
	\end{proof}
	
	\begin{theorem} \label{thm:rapid_variation_WAP}
		Suppose $f$ satisfies the hypotheses of Lemma~\ref{lem:elasticity}, and set
		$b(t) = F^{-1}(t)$. Then $t b'(t) = o(b(t))$ as $t\to\infty$ (equivalently,
		$F(x) = o(x/f(x))$ as $x\to 0^+$).
	\end{theorem}
	
	\begin{proof}
		Setting $x = b(t)$ (so $t = F(x)$ and $b'(t)=-f(x)$), we have
		\[
		\lim_{t \to \infty} \frac{t b'(t)}{b(t)}
		= -\lim_{x \to 0^+} \frac{F(x)}{x/f(x)}.
		\]
		Since $F(x)\to\infty$ and $x/f(x)\to\infty$ as $x\to 0^+$, l'H\^opital's rule
		gives
		\[
		\lim_{x \to 0^+} \frac{F(x)}{x/f(x)}
		= \lim_{x \to 0^+} \frac{-1/f(x)}{(f(x) - xf'(x))/(f(x))^2}
		= \lim_{x \to 0^+} \frac{-1}{1 - xf'(x)/f(x)} = 0,
		\]
		where the last equality uses Lemma~\ref{lem:elasticity}.
	\end{proof}
	
	Theorem~\ref{thm:rapid_variation_WAP} explains why the decoupled finite-window
	condition of Theorem~\ref{thm.22} is indispensable in the super-flat regime:
	the elasticity of $b(t)=F^{-1}(t)$ vanishes at infinity, so \eqref{eq.8} fails.
	
	\subsection{Monotonicity and log-concavity criteria}
	
	For nonlinearities growing like $x^2$ or slower near zero, it is natural to
	assume that $f$, $f'$, and $f''$ are increasing. A Taylor--Landau argument
	immediately yields a uniform lower bound on $Q$.
	
	\begin{theorem} \label{thm:Q_lower_half}
		Suppose $f(x) > 0$, $f'(x) > 0$, and $f''(x) > 0$ on $(0,\delta)$, with $f'$
		and $f''$ increasing and $f(0) = f'(0^+) = 0$. Then $Q(x)\ge 1/2$ on
		$(0,\delta)$.
	\end{theorem}
	
	\begin{proof}
		Fix $x\in(0,\delta)$ and let $y\in[0,x]$. By Taylor's theorem with Lagrange
		remainder and monotonicity of $f''$, there exists $\xi\in(x-y,x)$ such that
		\[
		0\le f(x - y) = f(x) - y f'(x) + \frac{y^2}{2} f''(\xi)
		\le f(x) - y f'(x) + \frac{y^2}{2} f''(x) =: P(y).
		\]
		The quadratic $P(y)$ attains its global minimum at $y_v = f'(x)/f''(x) > 0$.
		Since $f'(0^+)=0$ and $f''$ is increasing,
		\[
		f'(x) = \int_0^x f''(t)\,dt \le x f''(x),
		\]
		so $y_v \in [0,x]$. Evaluating $P(y_v)\ge 0$ gives
		\[
		0 \le f(x) - \frac{(f'(x))^2}{2f''(x)},
		\]
		whence $Q(x)\ge 1/2$.
	\end{proof}
	
	Even without log-concavity, $Q(x)$ cannot stay strictly above $1$ near zero.
	
	\begin{theorem} \label{thm:Q_liminf_one}
		Suppose $f(x) > 0$, $f'(x) > 0$, and $f''(x) > 0$ on $(0,\delta)$, with
		$f(0) = f'(0^+) = 0$. Then $\liminf_{x\to 0^+} Q(x)\le 1$.
	\end{theorem}
	
	\begin{proof}
		Suppose to the contrary that $q^*:=\liminf_{x\to 0^+} Q(x) > 1$. Recall that
		$k(x):=f(x)/f'(x)>0$ satisfies $k'(x)=1-Q(x)$. Choose $\epsilon\in(0,1)$ with
		$1-q^*+\epsilon<0$; then there exists $x_1\in(0,\delta)$ such that
		$k'(x)\le 1-q^*+\epsilon<0$ on $(0,x_1)$. Since $f'$ is increasing and
		$f(0)=0$, the mean value theorem gives $f(x)=xf'(\xi)\le xf'(x)$, so
		$0<k(x)\le x$ and hence $k(0^+)=0$. For $0<x<y<x_1$, the mean value theorem
		yields $k(y)-k(x)\le (1-q^*+\epsilon)(y-x)$, i.e.,
		$k(x)-(1-q^*+\epsilon)x \ge k(y)-(1-q^*+\epsilon)y$. Letting $x\to 0^+$ gives
		$0 \ge k(y)-(1-q^*+\epsilon)y$, contradicting $k(y)>0$ and $1-q^*+\epsilon<0$.
	\end{proof}
	
	A uniform upper bound $Q(x)\le 1$ is immediately provided by log-concavity.
	
	\begin{theorem} \label{thm:logconcave_Q}
		Suppose $f \in C^2((0,\delta))$ is log-concave on $(0,\delta)$ with $f(x)>0$
		and $f'(x)>0$. Then $Q(x)\le 1$ on $(0,\delta)$, and hence
		$\limsup_{x\to 0^+} Q(x)\le 1$.
	\end{theorem}
	
	\begin{proof}
		Log-concavity means $\frac{d^2}{dx^2}\log f(x) = \frac{f''(x)f(x)-(f'(x))^2}{(f(x))^2} \le 0$
		on $(0,\delta)$, which rearranges immediately to $Q(x)\le 1$.
	\end{proof}
	
	Combining Theorems~\ref{thm:Q_implies_WAP}, \ref{thm:Q_lower_half}, and
	\ref{thm:logconcave_Q}, we see that any $C^2$ log-concave function with
	$f(0)=f'(0^+)=0$ and increasing $f''>0$ automatically satisfies
	$\frac{1}{2}\le Q(x)\le 1$, and therefore yields a WAP composition
	$f\circ F^{-1}$.
	
	\subsubsection{Decreasing second derivative and elasticity above one}
	
	When $f''$ is decreasing near zero (typical of power-law growth $x^\beta$ with
	$\beta\in(1,2)$), $Q$ is bounded above by $1/2$ automatically, and bounded
	below whenever the elasticity $E(x):=xf'(x)/f(x)$ stays strictly above~$1$.
	
	\begin{theorem} \label{thm:Q_upper_half}
		Suppose $f$ and $f'$ are increasing on $(0,\delta)$, $f''$ is positive and
		non-increasing on $(0,\delta)$, and $f(0)=f'(0^+)=0$. Then $Q(x)\le 1/2$ on
		$(0,\delta)$.
	\end{theorem}
	
	\begin{proof}
		Fix $x\in(0,\delta)$ and put $w(s)=f''(s)\ge 0$. Since $f(0)=f'(0^+)=0$,
		\[
		f'(x)=\int_0^x w(u)\,du, \qquad f(x)=\int_0^x (x-s)w(s)\,ds.
		\]
		Because $w$ is non-increasing, $\int_s^x w(u)\,du \ge (x-s)w(x)$ for every
		$s\in[0,x]$. Multiplying by $w(s)$ and integrating in $s$ yields
		\[
		\int_0^x w(s)\int_s^x w(u)\,du\,ds
		\ge w(x)\int_0^x (x-s)w(s)\,ds = w(x)f(x).
		\]
		By symmetry on $[0,x]^2$, the left-hand side is
		\[
		\frac12\left(\int_0^x w(u)\,du\right)^2
		=\frac12(f'(x))^2.
		\]
		It follows that $Q(x)=f(x)f''(x)/(f'(x))^2\le 1/2$.
	\end{proof}
	
	\begin{theorem} \label{thm:Q_lower_elasticity}
		Let $f:[0,\delta)\to\mathbb{R}$ be twice differentiable on $(0,\delta)$ with
		$f(0)=f'(0^+)=0$, such that $f$ and $f'$ are strictly increasing and $f''$ is
		strictly positive and decreasing on $(0,\delta)$. Put $E(x):=xf'(x)/f(x)$. If
		$\liminf_{x\to 0^+} E(x) > 1$, then $\liminf_{x\to 0^+} Q(x) > 0$.
	\end{theorem}
	
	Note that the threshold $\liminf_{x\to 0^+} E(x) > 1$ is sharp: for
	$f(x)=x/\log(1/x)$, all other hypotheses hold with $E(x)\to 1$, yet
	$Q(x)\to 0$ and $f\circ F^{-1}$ is rapidly varying (non-WAP).
	
	\begin{proof}
		Let $p := \liminf_{x\to 0^+} E(x) > 1$, and set $A(x):=xf''(x)/f'(x)$, so that
		$Q(x)=A(x)/E(x)$. Since $f''$ is decreasing, $f'$ is concave with $f'(0^+)=0$,
		whence $f(x)=\int_0^x f'(t)\,dt \ge \frac{1}{2}xf'(x)$, i.e., $E(x)\le 2$ and
		$Q(x)\ge A(x)/2$.
		
		It therefore suffices to show $\liminf_{x\to 0^+} A(x)>0$. Suppose for
		contradiction that there exists a sequence $x_n\to 0^+$ with
		$A(x_n)=:\epsilon_n\to 0$, so $x_n f''(x_n)=\epsilon_n f'(x_n)$. Fix $K>1$ (to
		be chosen shortly) and set $y_n := K x_n$. Since $f''$ is decreasing,
		\[
		f'(y_n) = f'(x_n) + \int_{x_n}^{Kx_n} f''(t)\,dt
		\le f'(x_n) + (K-1)x_n f''(x_n) = f'(x_n)\big[1+\epsilon_n(K-1)\big],
		\]
		while $f'$ increasing gives
		\[
		f(y_n) = f(x_n) + \int_{x_n}^{Kx_n} f'(t)\,dt \ge (K-1)x_n f'(x_n).
		\]
		Consequently,
		\[
		E(y_n) = \frac{Kx_n f'(Kx_n)}{f(Kx_n)} \le \frac{K}{K-1}\big[1+\epsilon_n(K-1)\big],
		\]
		so $\limsup_{n\to\infty} E(y_n) \le \frac{K}{K-1}$. Choosing
		$K=(p+1)/(p-1)$ gives $K>1$ and $\frac{K}{K-1}<p$, contradicting
		$\liminf_{n\to\infty} E(y_n)\ge p$. Hence $\liminf_{x\to 0^+} A(x)=c>0$, and
		therefore $\liminf_{x\to 0^+} Q(x)\ge c/2>0$.
	\end{proof}
	
	\subsubsection{Why monotonicity alone is insufficient}
	
	In Theorem~\ref{thm:logconcave_Q}, log-concavity supplied the upper bound
	$Q(x)\le 1$. One might ask whether monotonicity of $f$, $f'$, and $f''$ (even
	strengthened to power-weighted monotonicity) already forces $Q$ to be bounded
	above. The following construction shows that it does not.
	
	\begin{theorem} \label{thm:Q_unbounded_monotone}
		There exists a smooth, strictly positive function $f: (0, 1] \to (0,\infty)$
		with $f'(x)>0$ and $f''(x)>0$ such that:
		\begin{enumerate}
			\item $f''$ is strictly increasing on $(0, 1]$;
			\item $x\mapsto f'(x)/x$ is strictly increasing on $(0, 1]$;
			\item $x\mapsto f(x)/x^{3/2}$ is strictly increasing on $(0, 1]$;
		\end{enumerate}
		while $\limsup_{x \to 0^+} Q(x) = +\infty$.
	\end{theorem}
	
	\begin{proof}
		We construct a piecewise-constant non-decreasing template $g_2$ for $f''$ on
		$(0,x_2]$, evaluate its primitives $g_1(x)=\int_0^x g_2(t)\,dt$ and
		$g_0(x)=\int_0^x (x-t)g_2(t)\,dt$, and then smooth $g_2$.
		
		\paragraph{Step 1: Structural setup.}
		Choose a sequence $\{x_k\}_{k=2}^\infty$ in $(0,1)$ decaying rapidly to $0$ so
		that
		\begin{equation}\label{eq:sep_cond}
			x_{k+1} \le \frac{x_k}{2k^2} \qquad \text{for all } k \ge 2.
		\end{equation}
		On each interval $(x_{k+1}, x_k]$, set $\delta_k := x_k/k$. By
		\eqref{eq:sep_cond}, $\delta_k < x_k - x_{k+1}$ for all $k\ge 2$. Define the
		piecewise-constant function $g_2$ on $(x_{k+1}, x_k]$ by
		\[
		g_2(t) = \begin{cases}
			m_k & \text{for } t \in (x_{k+1}, x_k - \delta_k], \\
			M_k := k\, m_k & \text{for } t \in (x_k - \delta_k, x_k],
		\end{cases}
		\]
		where the levels $\{m_k\}_{k=2}^\infty$ are prescribed recursively by
		\begin{equation}\label{eq:mk_rec}
			m_2 = 1, \qquad m_{k+1} = \frac{m_k}{k+1} \implies M_{k+1} = (k+1)m_{k+1} = m_k.
		\end{equation}
		Then $g_2$ is positive and non-decreasing on $(0, x_2]$, and for every $k\ge 2$,
		\begin{equation}\label{eq:hist_bound}
			g_2(t) \le M_{k+1} = m_k \qquad \text{for all } t \in (0, x_{k+1}].
		\end{equation}
		
		\paragraph{Step 2: Upper bound for $g_1(x_k)$.}
		Using \eqref{eq:hist_bound}, $M_k = k m_k$, and $\delta_k = x_k/k$, we have
		\begin{align*}
			g_1(x_k) = \int_0^{x_k} g_2(t)\,dt
			&= \int_0^{x_{k+1}} g_2(t)\,dt + \int_{x_{k+1}}^{x_k-\delta_k} m_k\,dt + \int_{x_k-\delta_k}^{x_k} M_k\,dt \\
			&\le m_k x_{k+1} + m_k(x_k-\delta_k-x_{k+1}) + k m_k \delta_k \\
			&= 2 m_k x_k - m_k \delta_k < 2 m_k x_k.
		\end{align*}
		Hence
		\begin{equation}\label{eq:g1_sq_bound}
			(g_1(x_k))^2 < 4 m_k^2 x_k^2.
		\end{equation}
		
		\paragraph{Step 3: Lower bound for $g_0(x_k)$.}
		Dropping the nonnegative contribution from $(0, x_{k+1}]$ in
		$g_0(x_k) = \int_0^{x_k}(x_k-t)g_2(t)\,dt$ gives
		\begin{align*}
			g_0(x_k)
			&\ge \int_{x_{k+1}}^{x_k-\delta_k} (x_k-t)m_k\,dt + \int_{x_k-\delta_k}^{x_k} (x_k-t)M_k\,dt \\
			&= \frac{m_k}{2}\big[(x_k-x_{k+1})^2-\delta_k^2\big] + \frac{M_k \delta_k^2}{2} \\
			&\ge \frac{m_k x_k^2}{2}\left[1 - \frac{2x_{k+1}}{x_k} - \frac{1}{k^2} + \frac{1}{k}\right].
		\end{align*}
		By \eqref{eq:sep_cond}, $2x_{k+1}/x_k \le 1/k^2$, so for all $k\ge 2$ the
		bracketed factor is at least $1 + \frac{1}{k} - \frac{2}{k^2} \ge 1$. Thus
		\begin{equation}\label{eq:g0_lower}
			g_0(x_k) \ge \frac{m_k x_k^2}{2}.
		\end{equation}
		
		\paragraph{Step 4: Divergence of the template quotient.}
		Since $g_2(x_k)=M_k=k m_k$, combining \eqref{eq:g1_sq_bound} and
		\eqref{eq:g0_lower} gives
		\[
		\frac{g_0(x_k) g_2(x_k)}{(g_1(x_k))^2}
		\ge \frac{\left(\frac{1}{2}m_k x_k^2\right)(k m_k)}{4 m_k^2 x_k^2} = \frac{k}{8} \longrightarrow \infty.
		\]
		
		\paragraph{Step 5: Smoothing and verification.}
		The jumps of $g_2$ occur at $a_k := x_k - \delta_k$ with heights
		$\Delta_k := M_k - m_k = (k-1)m_k$. Fix a non-decreasing $\vartheta\in C^\infty(\mathbb{R})$
		with $\vartheta(u)=0$ for $u\le -1$ and $\vartheta(u)=1$ for $u\ge 1$. Given
		$\varepsilon>0$, choose $\eta_k\in(0,\delta_k/2)$ so small that the intervals
		$(a_k-\eta_k, a_k+\eta_k)$ are pairwise disjoint and
		$2\eta_k\Delta_k \le \varepsilon 2^{-k} m_k x_k$. Replacing each jump at $a_k$
		by $m_k + \Delta_k \vartheta((x-a_k)/\eta_k)$ and extending by the constant
		$M_2$ on $[x_2,1]$ yields a non-decreasing function $\widetilde{g}_2\in C^\infty((0,1])$
		with $\widetilde{g}_2(x_k)=M_k$ and
		\[
		\int_0^{x_k} |\widetilde{g}_2(s)-g_2(s)|\,ds
		\le \sum_{j\ge k} 2\eta_j\Delta_j
		\le \varepsilon \sum_{j\ge k} 2^{-j} m_j x_j
		\le \varepsilon m_k x_k.
		\]
		Define
		\[
		f''(x) := \widetilde{g}_2(x) + \varepsilon x, \qquad
		f'(x) := \int_0^x f''(s)\,ds, \qquad
		f(x) := \int_0^x (x-s)f''(s)\,ds.
		\]
		From \eqref{eq:sep_cond} and \eqref{eq:mk_rec}, $x_k/m_k \le x_2 < 1$ for all
		$k\ge 2$. Hence the drift $\varepsilon x$ contributes at most
		$\varepsilon x_k^2/2 \le \varepsilon m_k x_k$ to $f'(x_k)$ and at most
		$\varepsilon x_k^3/6 \le \varepsilon m_k x_k^2/6$ to $f(x_k)$, while the
		smoothing error contributes at most $\varepsilon m_k x_k$ to $f'(x_k)$ and at
		most $\varepsilon m_k x_k^2$ to $f(x_k)$. Choosing $\varepsilon>0$ sufficiently
		small, there exist constants $c_1, c_2 > 0$ independent of $k$ such that
		\[
		f'(x_k) \le c_1 m_k x_k, \qquad f(x_k) \ge c_2 m_k x_k^2, \qquad f''(x_k) \ge k m_k,
		\]
		whence $Q(x_k) \ge (c_2/c_1^2)k \to \infty$.
		
		Finally, we verify the three monotonicity properties:
		\begin{enumerate}
			\item Since $\widetilde{g}_2$ is smoothly non-decreasing and $\varepsilon>0$,
			$f'''(x)=\widetilde{g}_2'(x)+\varepsilon>0$, so $f''$ is strictly increasing.
			\item Since $f''>0$ is strictly increasing, $f'(x)=\int_0^x f''(t)\,dt < x f''(x)$,
			so $\frac{d}{dx}\big(f'(x)/x\big) = (xf''(x)-f'(x))/x^2 > 0$.
			\item Since $f'(t)/t$ is strictly increasing, $f'(t) < t f'(x)/x$ for $0<t<x$,
			so $f(x)=\int_0^x f'(t)\,dt < \frac{1}{2}x f'(x)$. Hence
			\[
			\frac{d}{dx}\left(\frac{f(x)}{x^{3/2}}\right)
			= \frac{x f'(x) - \frac{3}{2}f(x)}{x^{5/2}}
			> \frac{\frac{1}{4}x f'(x)}{x^{5/2}} > 0.
			\]
		\end{enumerate}
		This completes the proof.
	\end{proof}
	
	\section{Auxiliary properties of WAP functions}
	\label{sec:properties_WAP}
	
	We now prove the two auxiliary lemmas for WAP functions stated in
	Section~\ref{sec:main_results}.
	
	\begin{proof}[Proof of Lemma~\ref{lemma.intwapiswap}]
		Since $\varphi\in L^1(0,\infty)$ is positive, $\phi(t):=\int_t^\infty \varphi(s)\,ds$
		is well defined and strictly decreasing. For $\lambda\in(0,1)$, monotonicity of
		the integral gives $\phi(\lambda t)\ge \phi(t)$, so
		$\liminf_{t\to\infty} \phi(\lambda t)/\phi(t) \ge 1 > 0$.
		
		For an upper bound valid for every $\lambda>0$, let
		$\overline{\Phi}_\lambda := \limsup_{t\to\infty} \varphi(\lambda t)/\varphi(t) < \infty$.
		For any $\varepsilon\in(0,1)$, there exists $T(\lambda,\varepsilon)>0$ such that
		$\varphi(\lambda s) < \overline{\Phi}_\lambda(1+\varepsilon)\varphi(s)$ for all
		$s\ge T(\lambda,\varepsilon)$. Substituting $s\mapsto \lambda s$ for
		$t\ge T(\lambda,\varepsilon)$ yields
		\[
		\frac{\phi(\lambda t)}{\phi(t)}
		= \frac{\lambda \int_t^\infty \varphi(\lambda s)\,ds}{\int_t^\infty \varphi(s)\,ds}
		\le \lambda \overline{\Phi}_\lambda (1+\varepsilon).
		\]
		Letting $t\to\infty$ and then $\varepsilon\to 0^+$ gives
		\[
		\Phi^*_\lambda := \limsup_{t\to\infty} \frac{\phi(\lambda t)}{\phi(t)}
		\le \lambda \overline{\Phi}_\lambda < +\infty \qquad \text{for all } \lambda>0.
		\]
		Finally, for $\lambda>1$, setting $\tau = \lambda t$ gives
		\[
		\liminf_{t\to\infty} \frac{\phi(\lambda t)}{\phi(t)}
		= \left(\limsup_{\tau\to\infty} \frac{\phi(\tau/\lambda)}{\phi(\tau)}\right)^{-1}
		= \frac{1}{\Phi^*_{1/\lambda}} \ge \frac{\lambda}{\overline{\Phi}_{1/\lambda}} > 0.
		\]
		Hence \eqref{eq.WAPlimitsInf} holds for all $\lambda>0$, so $\phi$ is WAP at
		$\infty$.
	\end{proof}
	
	\begin{proof}[Proof of Lemma~\ref{lemma.wapmeansnofasterthanpowerlaw}]
		Because $\varphi$ is WAP at $\infty$, $\liminf_{t\to\infty} \varphi(2t)/\varphi(t) > 0$,
		so there exist $c>0$ and $t_0>0$ such that $\varphi(2t)\ge c\varphi(t)$ for all
		$t\ge t_0$. Since $\varphi$ is decreasing, $c\in(0,1]$. Iterating gives
		$\varphi(2^k t_0)\ge c^k \varphi(t_0)$ for every integer $k\ge 0$.
		
		Given any $t>t_0$, let $k\ge 0$ be the unique integer with
		$2^k t_0 \le t < 2^{k+1} t_0$. Then $k\le \log(t/t_0)/\log 2$, and monotonicity
		of $\varphi$ gives
		\[
		\varphi(t) \ge \varphi(2^{k+1} t_0) \ge c^{k+1}\varphi(t_0)
		\ge c\cdot c^{\log(t/t_0)/\log 2}\varphi(t_0) = K t^\alpha,
		\]
		where $\alpha := \log c/\log 2 \le 0$ and $K := c t_0^{-\alpha}\varphi(t_0) > 0$.
		Therefore, for all $t\ge \max(e,t_0)$,
		\[
		\frac{\log \varphi(t)}{\log t} \ge \frac{\log K}{\log t} + \alpha,
		\]
		so $\lambda = \liminf_{t\to\infty} \frac{\log \varphi(t)}{\log t} \ge \alpha > -\infty$.
	\end{proof}
	
	\appendix
	
	\section{Sublinear $f\in\mathrm{RV}_0(1)$ with $f\circ F^{-1}$ non-WAP}
	\label{app:rv1}
	
	In this appendix we construct, for a broad class of sublinear functions
	$f\in\mathrm{RV}_0(1)$ (for which $f\circ F^{-1}$ is rapidly varying and hence
	non-WAP), an explicit positive integrable forcing $g$ whose solution preserves
	the rate $F(x(t))/t\to 1$ while violating both $g(t)=o((f\circ F^{-1})(t))$ and
	$\int_t^\infty g(s)\,ds = o(F^{-1}(t))$.
	
	\begin{hypothesis}[Standing hypotheses on $f\in\mathrm{RV}_0(1)$]\label{hyp:f_properties}
		Let $f:(0,1]\to(0,\infty)$ be $C^2$ and satisfy:
		\begin{enumerate}
			\item[\textbf{(H1)}] \textbf{Positivity and convexity:} $f(x)>0$, $f'(x)>0$, and $f''(x)>0$ for all $x\in(0,1]$.
			\item[\textbf{(H2)}] \textbf{Regular variation of index $1$:} $f\in\mathrm{RV}_0(1)$ with $\lim_{x\to 0^+} f(x)/x=0$ and $xf'(x)\sim f(x)$ as $x\to 0^+$.
			\item[\textbf{(H3)}] \textbf{Derivative regularity:} $xf''(x)=o(f'(x))$ as $x\to 0^+$.
		\end{enumerate}
	\end{hypothesis}
	
	Throughout this appendix, write
	\[
	h(t):=F^{-1}(t), \qquad m(t):=f'(h(t)).
	\]
	By Hypothesis~\ref{hyp:f_properties}, $m$ is positive, decreasing, and tends to
	zero, with $m'(t)=-f(h(t))f''(h(t))$. Define
	\[
	\delta(t):=-\frac{m'(t)}{(m(t))^2}=\frac{f(h(t))f''(h(t))}{(f'(h(t)))^2}.
	\]
	By~(H2) and~(H3),
	\begin{equation}\label{eq:delta_to_zero}
		\delta(t) = \frac{f(h(t))}{h(t)f'(h(t))}\cdot\frac{h(t)f''(h(t))}{f'(h(t))} \longrightarrow 0
		\qquad \text{as } t\to\infty.
	\end{equation}
	Moreover, substituting $x=h(s)$ (with $h'(s)=-f(h(s))$) shows that for every
	sufficiently large $T$,
	\begin{equation}\label{eq:m_integral_diverges}
		\int_T^t m(s)\,ds = \int_{h(t)}^{h(T)}\frac{f'(x)}{f(x)}\,dx = \log\frac{f(h(T))}{f(h(t))} \longrightarrow \infty
		\qquad \text{as } t\to\infty.
	\end{equation}
	
	\begin{theorem}[Construction of the self-regulating shift]\label{thm:direct_alpha_construction}
		Assume Hypothesis~\ref{hyp:f_properties}. Fix $T>0$ and choose
		$\alpha_T\in(0,T)$. There exists a unique global $C^1$ solution
		$\alpha:[T,\infty)\to(0,\infty)$ of
		\begin{equation}\label{eq:direct_alpha_ode}
			\alpha'(t)=\exp\!\left(-\frac12\alpha(t)m(t)\right),
			\qquad \alpha(T)=\alpha_T,
		\end{equation}
		possessing the following properties:
		\begin{enumerate}
			\item $0<\alpha'(t)<1$ and $\alpha(t)>0$ for all $t\geq T$;
			\item $\alpha(t)m(t)\to\infty$ as $t\to\infty$;
			\item $\alpha'(t)\to0$ and $\alpha(t)=o(t)$ as $t\to\infty$;
			\item $t-\alpha(t)\to\infty$ and $0<\alpha(t)<t$ for all $t\geq T$.
		\end{enumerate}
	\end{theorem}
	
	\begin{proof}
		The right-hand side of \eqref{eq:direct_alpha_ode} is continuous in $t$,
		locally Lipschitz in $\alpha$, positive, and bounded above by $1$. Hence a
		unique global $C^1$ solution exists on $[T,\infty)$ with $\alpha(t)>0$ and
		$0<\alpha'(t)<1$.
		
		Set $\Lambda(t):=\alpha(t)m(t)$. Differentiating and using
		\eqref{eq:direct_alpha_ode} together with $m'=-\delta m^2$ yields
		\begin{equation}\label{eq:L_differential_t}
			\Lambda'(t)=\alpha'(t)m(t)+\alpha(t)m'(t)=m(t)e^{-\Lambda(t)/2}-\delta(t)m(t)\Lambda(t).
		\end{equation}
		Introduce the strictly increasing time variable $u=u(t):=\int_T^t m(s)\,ds$,
		which tends to $\infty$ by \eqref{eq:m_integral_diverges}. Writing $\Lambda$
		and $\delta$ as functions of $u$, \eqref{eq:L_differential_t} becomes
		\begin{equation}\label{eq:L_differential_u}
			\frac{d\Lambda}{du}=e^{-\Lambda/2}-\delta(t(u))\Lambda.
		\end{equation}
		Fix any $K>0$. By \eqref{eq:delta_to_zero}, there exists $u_K$ such that
		$\delta(t(u))K \le \frac{1}{2}e^{-K/2}$ for all $u\ge u_K$. Whenever
		$u\ge u_K$ and $\Lambda(u)\le K$, \eqref{eq:L_differential_u} gives
		\[
		\frac{d\Lambda}{du} \ge e^{-K/2}-\delta(t(u))K \ge \frac{1}{2}e^{-K/2} > 0.
		\]
		Hence $\Lambda$ must reach $K$ in finite $u$-time and cannot thereafter cross
		below $K$. Since $K>0$ is arbitrary, $\Lambda(t)\to\infty$ as $t\to\infty$.
		
		From \eqref{eq:direct_alpha_ode}, $\alpha'(t)=e^{-\Lambda(t)/2}\to 0$, and
		averaging gives $\alpha(t)/t = \alpha(T)/t + \frac{1}{t}\int_T^t \alpha'(s)\,ds \to 0$.
		Finally, since $\alpha(T)<T$ and $\alpha'(t)<1$,
		\[
		t-\alpha(t) = T-\alpha(T) + \int_T^t (1-\alpha'(s))\,ds > 0 \qquad (t\ge T),
		\]
		and $t-\alpha(t)\sim t\to\infty$.
	\end{proof}
	
	\begin{theorem}[Rate-preserving trajectory with amplified perturbation tail]\label{thm:direct_trajectory_properties}
		Assume Hypothesis~\ref{hyp:f_properties}, and let $\alpha$ be given by
		Theorem~\ref{thm:direct_alpha_construction}. For $t\geq T$, define
		\begin{equation}\label{eq:direct_x_g}
			x(t):=F^{-1}(t-\alpha(t)), \qquad g(t):=\alpha'(t)f(F^{-1}(t-\alpha(t))).
		\end{equation}
		Then $x$ is a positive $C^1$ solution of $x'(t)=-f(x(t))+g(t)$ on $[T,\infty)$
		satisfying
		\begin{equation}\label{eq:direct_rate_preservation}
			\lim_{t\to\infty}\frac{F(x(t))}{t}=1.
		\end{equation}
		Moreover, $g$ is positive and integrable on $[T,\infty)$, and
		\begin{equation}\label{eq:direct_tail_amplification}
			\lim_{t\to\infty}\frac{g(t)}{f(F^{-1}(t))}=\infty, \qquad
			\lim_{t\to\infty}\frac{\int_t^\infty g(s)\,ds}{F^{-1}(t)}=\infty.
		\end{equation}
	\end{theorem}
	
	\begin{proof}
		Differentiating $F(x(t))=t-\alpha(t)$ gives $-x'(t)/f(x(t))=1-\alpha'(t)$, so
		$x'(t)=-f(x(t))+g(t)$. By Theorem~\ref{thm:direct_alpha_construction},
		$F(x(t))/t=(t-\alpha(t))/t\to 1$ and $g(t)>0$ on $[T,\infty)$.
		
		Since $\frac{d}{ds}\log f(F^{-1}(s)) = -f'(F^{-1}(s))=-m(s)$ and $m$ is
		decreasing, we obtain
		\begin{equation}\label{eq:direct_log_ratio}
			\log\frac{f(F^{-1}(t-\alpha(t)))}{f(F^{-1}(t))}
			= \int_{t-\alpha(t)}^t m(s)\,ds \ge \alpha(t)m(t) = \Lambda(t).
		\end{equation}
		Combining \eqref{eq:direct_alpha_ode} and \eqref{eq:direct_log_ratio} yields
		\begin{equation}\label{eq:direct_pointwise_amplification}
			\frac{g(t)}{f(F^{-1}(t))}
			= \alpha'(t)\frac{f(F^{-1}(t-\alpha(t)))}{f(F^{-1}(t))}
			\ge e^{-\Lambda(t)/2}e^{\Lambda(t)} = e^{\Lambda(t)/2} \longrightarrow \infty.
		\end{equation}
		
		To verify integrability of $g$, put $r(s):=s-\alpha(s)$. Since $\alpha'(s)\to 0$,
		we have $\alpha'(s)\le 1/2$ for $s\ge T_1$. Substituting $v=r(s)$ gives
		\[
		\int_{T_1}^\infty g(s)\,ds
		= \int_{r(T_1)}^\infty \frac{\alpha'(r^{-1}(v))}{1-\alpha'(r^{-1}(v))} f(F^{-1}(v))\,dv
		\le \int_{r(T_1)}^\infty f(F^{-1}(v))\,dv = F^{-1}(r(T_1)) < \infty.
		\]
		Finally, since $\int_t^\infty f(F^{-1}(s))\,ds = F^{-1}(t)\to 0$ and
		$\int_t^\infty g(s)\,ds\to 0$, l'H\^opital's rule together with
		\eqref{eq:direct_pointwise_amplification} gives
		\[
		\lim_{t\to\infty}\frac{\int_t^\infty g(s)\,ds}{F^{-1}(t)}
		= \lim_{t\to\infty}\frac{g(t)}{f(F^{-1}(t))} = \infty. \qedhere
		\]
	\end{proof}
	
	\begin{remark}
		The construction is self-regulating: $\Lambda(t)=\alpha(t)f'(F^{-1}(t))$ is
		driven to $\infty$ by $\int m = \infty$ and $-m'/m^2\to 0$, while choosing
		$\alpha'=e^{-\Lambda/2}$ forces $\alpha'\to 0$ slowly enough that the
		amplification $e^{\Lambda}$ in $f(F^{-1}(t-\alpha(t)))/f(F^{-1}(t))$ dominates
		$\alpha'$.
	\end{remark}
	
	\section{Asymptotically linear $f$ with $f\circ F^{-1}$ non-WAP}
	\label{app:linear}
	
	In this appendix we treat the asymptotically linear regime $f(x)\sim ax$ as
	$x\to 0^+$ ($a>0$). We establish the exact rate-preservation characterisation
	$\liminf_{t\to\infty} F(G(t))/t \ge 1$ (Theorems~\ref{thm:linear_suff}
	and~\ref{thm:linear_nec}) and construct a rate-preserving perturbation with
	$G(t)/F^{-1}(t)\to\infty$ (Theorem~\ref{thm:linear_supercrit}).
	
	\subsection{Sufficiency}
	
	\begin{theorem}\label{thm:linear_suff}
		Let $x' = -f(x) + g(t)$ for $t \ge 0$ with $x(0) = \xi > 0$. Assume
		$f\in C([0,\infty);[0,\infty))$ with $f(x) > 0$ for $x > 0$ and $f(x) \sim ax$
		as $x \to 0^+$ ($a > 0$). Assume $g(t) > 0$, $g \in L^1(0, \infty)$, and
		\[
		\limsup_{t\to\infty} \frac{\log G(t)}{t} \le -a,
		\]
		where $G(t) = \int_t^\infty g(s)\,ds$. Then
		\[
		\lim_{t\to\infty} \frac{\log x(t)}{t} = -a.
		\]
	\end{theorem}
	
	\begin{proof}
		\textbf{Step 1: Asymptotic stability ($x(t) \to 0$).}
		Since $x(0)>0$ and $g(t)>0$, $x$ cannot cross zero (at any first zero,
		$x'=g>0$), so $x(t)>0$ and $f(x(t))>0$ for all $t\ge 0$. Set
		$H(t):=x(t)-\int_0^t g(s)\,ds$. Then $H'(t)=-f(x(t))<0$, so $H$ is decreasing,
		and $H(t)>-\|g\|_{L^1}>-\infty$. Hence $H(t)$ converges to a finite limit $L$,
		and therefore $x(t)\to L+\|g\|_{L^1}=:x_\infty\ge 0$. If $x_\infty>0$, then by
		continuity and positivity of $f$ there exist $c>0$ and $T>0$ such that
		$f(x(s))\ge c$ for $s\ge T$; integrating \eqref{eq:main_ode} over $[T,t]$ gives
		$x(t)-x(T)\le -c(t-T)+\|g\|_{L^1(T,\infty)}\to -\infty$, contradicting
		$x(t)>0$. Thus $x_\infty=0$.
		
		\textbf{Step 2: Lower bound ($\liminf$).}
		Given $\epsilon>0$, choose $t_1>0$ such that $f(x(t))\le(a+\epsilon)x(t)$ for
		$t\ge t_1$. Since $g>0$, $x'(t)>-(a+\epsilon)x(t)$, so
		$x(t)>x(t_1)e^{-(a+\epsilon)(t-t_1)}$ for $t\ge t_1$. Taking logarithms,
		dividing by $t$, letting $t\to\infty$, and then $\epsilon\to 0^+$ yields
		$\liminf_{t\to\infty} \frac{\log x(t)}{t}\ge -a$.
		
		\textbf{Step 3: Upper bound ($\limsup$).}
		Given $\epsilon\in(0,a)$, set $b:=a-\epsilon>0$ and choose $t_2\ge 0$ such that
		$f(x(t))\ge bx(t)$ for $t\ge t_2$. Then $x'(t)\le -bx(t)+g(t)$, so variation of
		constants gives
		\[
		x(t) \le x(t_2)e^{-b(t-t_2)} + \int_{t_2}^t e^{-b(t-s)}g(s)\,ds, \qquad t\ge t_2.
		\]
		Integrating by parts with $g(s)=-G'(s)$ yields
		\[
		\int_{t_2}^t e^{-b(t-s)}g(s)\,ds
		= e^{-b(t-t_2)}G(t_2) - G(t) + b\int_{t_2}^t e^{-b(t-s)}G(s)\,ds.
		\]
		By hypothesis, taking $\delta=\epsilon/2$ (so $a-\delta=b+\epsilon/2$), there
		exist $T\ge t_2$ and $C>0$ such that $G(s)\le Ce^{-(b+\epsilon/2)s}$ for
		$s\ge T$. Splitting the integral at $T$,
		\[
		b\int_{t_2}^t e^{-b(t-s)}G(s)\,ds
		\le b e^{-bt}\left[\int_{t_2}^T e^{bs}G(s)\,ds + C\int_T^t e^{-(\epsilon/2)s}\,ds\right]
		\le M e^{-bt},
		\]
		for a constant $M<\infty$. Dropping $-G(t)<0$ gives $x(t)\le \tilde{C}e^{-bt}$
		for $t\ge T$, whence $\limsup_{t\to\infty}\frac{\log x(t)}{t}\le -(a-\epsilon)$.
		Letting $\epsilon\to 0^+$ yields $\limsup_{t\to\infty}\frac{\log x(t)}{t}\le -a$,
		and combining Steps~2 and~3 completes the proof.
	\end{proof}
	
	\subsection{Necessity}
	
	\begin{theorem}\label{thm:linear_nec}
		Let $x' = -f(x) + g(t)$ for $t \ge 0$ with $x(0) = \xi > 0$, $g(t) > 0$, and
		$f(x) \sim ax$ as $x \to 0^+$ ($a > 0$). If
		\[
		\lim_{t\to\infty} \frac{\log x(t)}{t} = -a,
		\]
		then $g \in L^1(0, \infty)$ and
		\[
		\limsup_{t\to\infty} \frac{\log G(t)}{t} \le -a,
		\]
		where $G(t) = \int_t^\infty g(s)\,ds$.
	\end{theorem}
	
	\begin{proof}
		Fix $\epsilon\in(0,a)$. By hypothesis, there exists $T_1>0$ such that
		$x(t)\le e^{-(a-\epsilon)t}$ for $t\ge T_1$, and since $f(x)\sim ax$, there
		exists $T_2\ge T_1$ such that $f(x(t))\le (a+\epsilon)x(t)\le (a+\epsilon)e^{-(a-\epsilon)t}$
		for $t\ge T_2$. Thus $f\circ x \in L^1(0,\infty)$. Integrating
		\eqref{eq:main_ode} over $[0,t]$ gives
		$\int_0^t g(s)\,ds = x(t)-x(0)+\int_0^t f(x(s))\,ds$; letting $t\to\infty$
		shows that $g\in L^1(0,\infty)$.
		
		Integrating \eqref{eq:main_ode} over $[t,\infty)$ and using $x(s)\to 0$ yields
		\[
		G(t) = x(t) + \int_t^\infty f(x(s))\,ds.
		\]
		For $t\ge T_2$,
		\[
		G(t) \le e^{-(a-\epsilon)t} + (a+\epsilon)\int_t^\infty e^{-(a-\epsilon)s}\,ds
		= \left(1 + \frac{a+\epsilon}{a-\epsilon}\right)e^{-(a-\epsilon)t}.
		\]
		Taking logarithms, dividing by $t$, letting $t\to\infty$, and then
		$\epsilon\to 0^+$ gives $\limsup_{t\to\infty} \frac{\log G(t)}{t}\le -a$.
	\end{proof}
	
	\subsection{Supercritical perturbation}
	
	\begin{proposition}\label{prop:linear_k}
		Let $j: [t_0, \infty) \to \mathbb{R}$ be an increasing $C^1$ function such
		that $\lim_{t \to \infty} j(t)/t = 1$. Then there exists a $C^1$ function
		$k: [t_0, \infty) \to \mathbb{R}$ such that $k(t) \to \infty$, $k(t) = o(t)$
		as $t \to \infty$, and $j(t) - k(t)$ is non-decreasing on $[t_0, \infty)$.
	\end{proposition}
	
	\begin{proof}
		Define $h(u):=(1+\sqrt{u-j(t_0)+1})^{-1}\in(0,1/2]$ on $[j(t_0),\infty)$, so
		$h(u)\to 0$ and $\int_{j(t_0)}^\infty h(v)\,dv=\infty$. Set
		$H(u):=\int_{j(t_0)}^u h(v)\,dv$ and $k(t):=H(j(t))$. Then $k\in C^1$,
		$(j-k)'(t)=j'(t)(1-h(j(t)))\ge \frac{1}{2}j'(t)\ge 0$, $k(t)\to\infty$, and
		l'H\^opital's rule gives $H(u)/u\to 0$, whence
		$k(t)/t=(H(j(t))/j(t))(j(t)/t)\to 0$.
	\end{proof}
	
	\begin{lemma}\label{lem:linear_G}
		Let $a > 0$ and $f\in C((0,1];(0,\infty))$ satisfy $f(x)\sim ax$ as $x\to 0^+$.
		Set $j(t):=\frac{1}{a}\log(1/F^{-1}(t))$ and $G(t):=e^{-a(j(t)-k(t))}$, where
		$k$ is given by Proposition~\ref{prop:linear_k}. Then:
		\begin{enumerate}
			\item[(a)] $j\in C^1$, $j$ is strictly increasing for large $t$, and $j(t)/t\to 1$ as $t\to\infty$;
			\item[(b)] $G$ is ultimately non-increasing;
			\item[(c)] $\lim_{t\to\infty} \frac{\log G(t)}{t} = -a$;
			\item[(d)] $\lim_{t\to\infty} \frac{G(t)}{F^{-1}(t)} = \infty$.
		\end{enumerate}
	\end{lemma}
	
	\begin{proof}
		Since $F'(x)=-1/f(x)$, differentiating gives
		$j'(t)=\frac{1}{a}\frac{f(F^{-1}(t))}{F^{-1}(t)}\to 1$ as $t\to\infty$, so
		$j(t)/t\to 1$ by l'H\^opital's rule and $j$ is strictly increasing for large
		$t$, proving~(a). Next, $G'(t)=-a(j'(t)-k'(t))e^{-a(j(t)-k(t))}\le 0$ for large
		$t$ because $j-k$ is non-decreasing, proving~(b). Part~(c) follows from
		$\frac{\log G(t)}{t}=-a(j(t)/t - k(t)/t)\to -a$, and~(d) follows from
		$G(t)/F^{-1}(t)=e^{-a(j(t)-k(t))}/e^{-aj(t)}=e^{ak(t)}\to\infty$.
	\end{proof}
	
	\begin{theorem}\label{thm:linear_supercrit}
		Suppose $f\in C([0,\infty);[0,\infty))$ with $f(x)>0$ for $x>0$ satisfies
		$f(x)\sim ax$ as $x\to 0^+$ ($a>0$). Then there exists an ultimately positive
		continuous forcing $g\in L^1(0,\infty)$ such that
		\[
		\lim_{t\to\infty} \frac{G(t)}{F^{-1}(t)}=+\infty, \qquad
		\liminf_{t\to\infty} \frac{F(G(t))}{t}\ge 1,
		\]
		and the solution $x$ of $x'(t)=-f(x(t))+g(t)$ satisfies
		$\lim_{t\to\infty} F(x(t))/t=1$.
	\end{theorem}
	
	\begin{proof}
		Let $G$ be given by Lemma~\ref{lem:linear_G} and set $g(t):=-G'(t)$ for large
		$t$ (extended continuously and positively to $[0,\infty)$). Note that
		$(j-k)'(t)\ge \frac{1}{2}j'(t)>0$ for large $t$, so $g$ is ultimately strictly
		positive, continuous, and $G(t)=\int_t^\infty g(s)\,ds$. Since
		$\log G(t)/t\to -a$, Theorem~\ref{thm:linear_suff} yields $\log x(t)/t\to -a$.
		Because $F(x)\sim \frac{1}{a}\log(1/x)$ as $x\to 0^+$, it follows that
		$F(x(t))/t\to 1$ and $F(G(t))/t\to 1$, while $G(t)/F^{-1}(t)\to\infty$ by
		Lemma~\ref{lem:linear_G}(d).
	\end{proof}
	
	\begin{remark}\label{rem:linear_nonwap}
		When $f(x)\sim ax$ ($a>0$), $F(x)\sim \frac{1}{a}\log(1/x)$ is slowly varying
		at $0^+$, so $F^{-1}$ and $f\circ F^{-1}$ are rapidly varying at $\infty$ and
		hence non-WAP.
	\end{remark}
	
	\section{Superlinear $f$ with $f\circ F^{-1}$ non-WAP}
	\label{app:superlinear}
	
	We show here that if $f$ is superlinear at zero ($f(x)/x\to\infty$ as
	$x\to 0^+$) with $F(x)\to\infty$, then $f\circ F^{-1}$ is necessarily rapidly
	varying (hence non-WAP), and rate-preserving perturbations exist for which
	$G(t)/F^{-1}(t)\to\infty$.
	
	\begin{theorem}\label{thm:superlinear_nonwap}
		Let $f:[0,\infty)\to[0,\infty)$ be continuous, strictly increasing, with
		$f(0)=0$ and $f(x)>0$ for $x>0$. Suppose that
		\[
		\lim_{x\to 0^+} \frac{f(x)}{x}=\infty
		\]
		and $F(x)\to\infty$ as $x\to 0^+$. Then $\varphi:=f\circ F^{-1}$ is rapidly
		varying at $\infty$, and in particular is not WAP at $\infty$.
	\end{theorem}
	
	\begin{proof}
		Since $F$ is a strictly decreasing $C^1$ bijection from $(0,1]$ onto
		$[0,\infty)$, its inverse $F^{-1}$ is strictly decreasing with
		$F^{-1}(t)\to 0^+$, and $\varphi = f\circ F^{-1}$ is strictly decreasing. Fix
		$\lambda>1$ and for $t>0$ put $x_1=F^{-1}(t)$ and $x_2=F^{-1}(\lambda t)$, so
		$0<x_2<x_1\le 1$. Because $f$ is increasing,
		\[
		(\lambda-1)t = F(x_2)-F(x_1) = \int_{x_2}^{x_1}\frac{du}{f(u)}
		\le \frac{x_1-x_2}{f(x_2)} < \frac{x_1}{f(x_2)}.
		\]
		Rearranging and dividing by $\varphi(t)=f(x_1)$ gives
		\[
		0 < \frac{\varphi(\lambda t)}{\varphi(t)} = \frac{f(x_2)}{f(x_1)}
		< \frac{x_1}{f(x_1)}\cdot\frac{1}{(\lambda-1)t}.
		\]
		As $t\to\infty$, $x_1\to 0^+$, so $x_1/f(x_1)\to 0$ and $1/((\lambda-1)t)\to 0$.
		Hence $\varphi(\lambda t)/\varphi(t)\to 0$ for every $\lambda>1$: $\varphi$ is
		rapidly varying at $\infty$ and violates \eqref{eq.WAPlimitsInf}.
	\end{proof}
	
	\begin{theorem}\label{thm:superlinear_counterex}
		Let $f:[0,\infty)\to[0,\infty)$ be continuous, differentiable on $(0,\infty)$,
		and strictly increasing with $f(0)=0$, $f(x)>0$ for $x>0$, $F(x)\to\infty$ as
		$x\to 0^+$, and $\lim_{x\to 0^+} f'(x)=+\infty$. Then there exists a continuous,
		positive, integrable function $g$ such that
		\[
		\lim_{t\to\infty}\frac{\int_t^\infty g(s)\,ds}{F^{-1}(t)}=\infty,
		\]
		while the solution of $x'(t)=-f(x(t))+g(t)$ ($t\ge 0$) satisfies
		$F(x(t))/t\to 1$ and $\liminf_{t\to\infty} F(G(t))/t \ge 1$.
	\end{theorem}
	
	\begin{proof}
		Define $\alpha(t):=(1+t)^{1/2}-1$ for $t\ge 0$, $x(t):=F^{-1}(t-\alpha(t))$,
		and $g(t):=\alpha'(t)f(x(t))$. Then $x$ is positive and $C^1$, $g$ is
		continuous and positive, $x'(t)=-f(x(t))+g(t)$, and
		$F(x(t))/t=(t-\alpha(t))/t\to 1$.
		
		Set $\phi(t):=1/(f\circ F^{-1})(t)$ and $\psi(t):=\log\phi(t)$, so
		$\psi'(t)=f'(F^{-1}(t))=:\mu(t)\to\infty$ as $t\to\infty$. Since
		$\alpha'(t)\le 1/2$ and $f$ is increasing,
		$g(t)\le \frac{1}{2}f(F^{-1}(t/2))=\frac{1}{2\phi(t/2)}=:\tilde{\phi}(t)/2$.
		Because $\tilde{\phi}'(t)/\tilde{\phi}(t)=-\mu(t/2)/2\to -\infty$,
		$\tilde{\phi}$ decays faster than every negative exponential and is integrable;
		hence $g\in L^1(0,\infty)$.
		
		Now, by the mean value theorem, there exists $t^*(t)\in[t-\alpha(t),t]$ such
		that
		\[
		\log\frac{\phi(t)}{\phi(t-\alpha(t))} = \psi(t)-\psi(t-\alpha(t)) = \mu(t^*(t))\alpha(t).
		\]
		Since $t-\alpha(t)\to\infty$ and $\mu\to\infty$, we have $\mu(t^*(t))>1$ for
		$t\ge T_1$, so $\phi(t)/\phi(t-\alpha(t))>e^{\alpha(t)}$. For $t\ge T_2$,
		$\alpha(t)>t^{1/2}/2$ and $\alpha'(t)\ge t^{-1/2}/4$, whence
		\[
		\frac{g(t)}{f(F^{-1}(t))} = \frac{\alpha'(t)\phi(t)}{\phi(t-\alpha(t))}
		> \frac{1}{4}t^{-1/2}e^{t^{1/2}/2} \longrightarrow \infty.
		\]
		By l'H\^opital's rule applied to $G(t)/F^{-1}(t) = \int_t^\infty g(s)\,ds / \int_t^\infty f(F^{-1}(s))\,ds$,
		we obtain $G(t)/F^{-1}(t)\to\infty$.
		
		Finally, for any $\epsilon\in(0,1/2)$, choose $T(\epsilon)>0$ such that
		$\alpha'(t)<\epsilon$ and $\alpha(t)<\epsilon t$ for $t\ge T(\epsilon)$. Then
		\[
		G(t) = \int_t^\infty \alpha'(s)f(F^{-1}(s-\alpha(s)))\,ds
		< \epsilon \int_t^\infty f(F^{-1}((1-\epsilon)s))\,ds
		= \frac{\epsilon}{1-\epsilon} F^{-1}((1-\epsilon)t).
		\]
		Setting $x=F^{-1}((1-\epsilon)t)$ (so $t=F(x)/(1-\epsilon)$) and using
		$\epsilon/(1-\epsilon)<1$ with $F$ decreasing gives
		\[
		\liminf_{t\to\infty}\frac{F(G(t))}{t}
		\ge (1-\epsilon)\liminf_{x\to 0^+}\frac{F(\epsilon x/(1-\epsilon))}{F(x)}
		\ge 1-\epsilon.
		\]
		Letting $\epsilon\to 0^+$ gives $\liminf_{t\to\infty} F(G(t))/t \ge 1$.
	\end{proof}
	
	\section{Super-flat nonlinearities: failure of the tail criterion and of integrability under $x\sim F^{-1}$}
	\label{app:superflat}
	
	This appendix establishes two sharpness results for super-flat nonlinearities
	whose transformed function $\varphi = f\circ F^{-1}$ \emph{is} WAP at $\infty$
	(in fact $\varphi\in\mathrm{RV}_\infty(-1)$), while the two-sided elasticity
	condition \eqref{eq.8} of Theorem~\ref{thm.23} fails. First, for
	$f(x)=e^{-1/x}$, we construct a $C^1$ stalling perturbation satisfying
	$G(t)=o(F^{-1}(t))$ that destroys rate preservation. Second, for
	$f(x)=x^2 e^{-1/x}$, we construct a non-integrable forcing with
	$x(t)\sim F^{-1}(t)$, justifying our choice of $F(x(t))/t\to 1$ in
	Section~\ref{sec:main_results}.
	
	\subsection{Explicit bounds and failure of $G=o(F^{-1})$ via stalling}
	
	Consider $f(x)=e^{-1/x}$ for $x>0$, with $F(x)=\int_x^1 e^{1/u}\,du$.
	
	\begin{lemma}[Fundamental bounds on $F$ and $F^{-1}$]\label{lem:explicit_bounds}
		For all $x\in(0,0.1]$,
		\[
		\frac12 x^2 e^{1/x} \le F(x) \le 3 x^2 e^{1/x}.
		\]
		Furthermore, for every $v\ge e^{10}$,
		\[
		\frac{1}{2\ln v} \le F^{-1}(v) \le \frac{4}{\ln v}.
		\]
	\end{lemma}
	
	\begin{proof}
		Put $V=1/x$ and $H(v)=e^v/v^2$. Then $F(x)=\int_1^V H(v)\,dv$ and
		$H'(v)=H(v)(1-2/v)$. For $V\ge 10$, $H$ is increasing on $[4,V]$ with
		$\frac{1}{2}H(v)\le H'(v)\le H(v)$, so
		\[
		F(x) \ge \int_4^V H'(v)\,dv = H(V)-H(4) \ge \frac12 H(V),
		\]
		and
		\[
		F(x) \le \int_1^4 H(v)\,dv + 2\int_4^V H'(v)\,dv \le C_0 + 2H(V) \le 3H(V),
		\]
		since $C_0=\int_1^4 H(v)\,dv < H(10) \le H(V)$.
		
		For the inverse estimates, let $v=F(x)\ge e^{10}$, so $V=1/x>10$. Because $H$
		is increasing on $[10,\infty)$ and $3H(12)=3e^{12}/144<e^{10}$, the bound
		$v\le 3H(V)$ forces $V>12$. Then $v\ge \frac{1}{2}e^V/V^2$ yields
		$\ln v \ge V - 2\ln V - \ln 2 \ge V/2$ (for $V\ge 12$), giving
		$x=1/V \ge 1/(2\ln v)$. Conversely, $v\le 3e^V/V^2\le 3e^V$ implies
		$V\ge \ln v - \ln 3 \ge \frac{1}{4}\ln v$ for $v\ge e^{10}$, giving
		$x\le 4/\ln v$.
	\end{proof}
	
	\begin{theorem}[Insufficiency of $G=o(F^{-1})$ for $f(x)=e^{-1/x}$]\label{thm:stalling_counterex}
		Let $f(x)=e^{-1/x}$ for $x>0$, and fix any $c\in(0,1)$. There exists a
		continuous nonnegative forcing $g\in L^1(0,\infty)$ and a positive $C^1$
		solution $x$ of $x'(t)=-f(x(t))+g(t)$ such that
		\[
		\lim_{t\to\infty}\frac{\int_t^\infty g(s)\,ds}{F^{-1}(t)}=0,
		\]
		while
		\[
		\liminf_{t\to\infty}\frac{F(x(t))}{t}=c < 1 = \limsup_{t\to\infty}\frac{F(x(t))}{t}.
		\]
	\end{theorem}
	
	\begin{proof}
		Define the $C^1$ mollifier $\phi(u):=3u^2-2u^3$ on $[0,1]$ (so $\phi(0)=0$,
		$\phi(1)=1$, $\phi'(0)=\phi'(1)=0$, and $\int_0^1 \phi(u)\,du = 1/2$). Set
		$t_k := \exp(2^k)$ and choose an integer $k_0\ge 10$ large enough that
		$t_{k_0} \ge \max(100/c^2, 2e^{10})$. We construct a $C^1$ clock function $H$
		on $[t_{k_0},\infty)$ with $H(t_{k_0})=t_{k_0}$ by prescribing its derivative
		$h(t)=H'(t)\in[0,1]$ inductively on each interval $[t_k, t_{k+1}]$ ($k\ge k_0$).
		
		Given $v_k := H(t_k)$ (with $v_{k_0}=t_{k_0}$), set the target stall endpoint
		$\tau_k := (v_k + 1/2)/c$. Note inductively that
		$v_k \in [\frac{1}{2}t_k, t_k]$, so $t_k+1 < \tau_k < \tau_k+1 < t_{k+1}$ for
		all $k\ge k_0$. Divide $[t_k, t_{k+1}]$ into four phases:
		\begin{enumerate}
			\item \emph{Entering stall} $[t_k, t_k+1]$: $h(t) = 1 - \phi(t-t_k)$;
			\item \emph{Deep stall} $[t_k+1, \tau_k]$: $h(t) = 0$;
			\item \emph{Exiting stall} $[\tau_k, \tau_k+1]$: $h(t) = \phi(t-\tau_k)$;
			\item \emph{Drift phase} $[\tau_k+1, t_{k+1}]$: $h(t) = 1$.
		\end{enumerate}
		Then $h$ is continuous on $[t_{k_0},\infty)$ with $h(t)\in[0,1]$ and
		$h(t_{k_0})=1$. Extending $H$ to $[0,t_{k_0}]$ by setting $H(t)=t$ (so
		$h(t)=1$ and $g(t)=0$ on $[0,t_{k_0}]$) yields a $C^1$ non-decreasing clock
		function on $[0,\infty)$. Define $x(t):=F^{-1}(H(t))$ and
		$g(t):=f(x(t))[1-h(t)]\ge 0$ for $t\ge 0$; differentiating confirms
		$x'=-f(x)+g$ on $[0,\infty)$.
		
		Put $\alpha(t):=H(t)/t$. During $[t_k, t_k+1]$, $H$ increases by
		$\int_0^1 (1-\phi(u))\,du = 1/2$, so $H(t)=v_k+1/2$ throughout the deep stall
		$[t_k+1,\tau_k]$. In particular, $\alpha$ decreases on $[t_k+1,\tau_k]$ and
		attains $\alpha(\tau_k) = (v_k+1/2)/\tau_k = c$. During $[\tau_k,\tau_k+1]$,
		$H$ increases by $1/2$, and during the drift phase $[\tau_k+1,t_{k+1}]$ it
		increases by $t_{k+1}-(\tau_k+1)$, so
		\[
		v_{k+1} = H(t_{k+1}) = t_{k+1} - (\tau_k - v_k) \ge t_{k+1} - \left(\frac{1}{c}+1\right)t_k.
		\]
		Hence
		\[
		1 \ge \alpha(t_{k+1}) \ge 1 - \left(\frac{1}{c}+1\right)\frac{t_k}{t_{k+1}}
		= 1 - \left(\frac{1}{c}+1\right)\exp(-2^k).
		\]
		This confirms $v_{k+1}\ge \frac{1}{2}t_{k+1}$ inductively, shows
		$\limsup_{t\to\infty}\alpha(t)=1$, and (since $\alpha(t)\ge c$ on
		$[\tau_k,\tau_k+1]$ up to $O(1/\tau_k)$ and increases on $[\tau_k+1,t_{k+1}]$)
		yields $\liminf_{t\to\infty}\alpha(t)=c$.
		
		Now let $\Delta_k := \int_{t_k}^{t_{k+1}} g(s)\,ds$. Since $g(s)=0$ on the
		drift phase, $\Delta_k = \int_{t_k}^{\tau_k+1} f(F^{-1}(H(s)))[1-h(s)]\,ds$.
		Because $H(s)\ge v_k \ge \frac{1}{2}t_k \ge e^{10}$ on $[t_k,\tau_k+1]$ and
		$f\circ F^{-1}$ is decreasing, Lemma~\ref{lem:explicit_bounds} applied to
		$x_k=F^{-1}(v_k)$ gives
		\[
		f(x_k) \le \frac{3 x_k^2}{v_k} \le \frac{48}{v_k(\ln v_k)^2}.
		\]
		Bounding the stall duration by $\tau_k+1-t_k \le 2v_k/c$ and noting that
		$\ln v_k \ge \ln(t_k/2) = 2^k - \ln 2 \ge 2^{k-1}$ gives
		\[
		\Delta_k \le \left(\frac{2v_k}{c}\right)\left(\frac{48}{v_k(\ln v_k)^2}\right)
		= \frac{96/c}{(\ln v_k)^2} \le \frac{384/c}{4^k} =: \frac{M}{4^k}.
		\]
		For any $t\ge t_{k_0}$, choose $m\ge k_0$ such that $t\in[t_m, t_{m+1})$. Then
		\[
		\int_t^\infty g(s)\,ds \le \sum_{k=m}^\infty \Delta_k \le \frac{4M}{3\cdot 4^m},
		\]
		whereas Lemma~\ref{lem:explicit_bounds} and $t < t_{m+1}=\exp(2^{m+1})$ yield
		$F^{-1}(t) \ge \frac{1}{2\ln t} \ge \frac{1}{4\cdot 2^m}$. Therefore,
		\[
		\frac{\int_t^\infty g(s)\,ds}{F^{-1}(t)} \le \frac{16M}{3\cdot 2^m} \longrightarrow 0
		\qquad \text{as } t\to\infty. \qedhere
		\]
	\end{proof}
	
	\subsection{Non-integrable forcing with $x(t)\sim F^{-1}(t)$}
	
	\begin{proposition}\label{prop:solution-asymptotic-not-tail}
		We exhibit an odd continuous function $f$, strictly increasing and positive
		on $(0,1]$, a continuous forcing $g$, strictly positive for all sufficiently
		large times, and a positive solution of $x'(t)=-f(x(t))+g(t)$ such that
		$\varphi = f\circ F^{-1}\in \mathrm{RV}_\infty(-1)$ is WAP at $\infty$,
		\[
		x(t)\sim F^{-1}(t) \qquad \text{as } t\to\infty,
		\]
		and yet $\int_t^\infty g(s)\,ds=\infty$ for every sufficiently large $t$.
	\end{proposition}
	
	\begin{proof}
		For $0<z\leq 1$, set $f(z):=z^2 e^{-1/z}$, extend $f$ oddly to $z<0$, and set
		$f(0)=0$. Then $f$ is continuous and strictly increasing on $(0,1]$, with
		escape-time function
		\[
		F(z)=\int_z^1\frac{dw}{w^2 e^{-1/w}}=e^{1/z}-e, \qquad 0<z\leq 1.
		\]
		Consequently,
		\[
		u(t)=F^{-1}(t)=\frac{1}{\log(t+e)}, \qquad
		\varphi(t)=f(u(t))=\frac{1}{(t+e)[\log(t+e)]^2},
		\]
		so $\varphi\in\mathrm{RV}_\infty(-1)$.
		
		For sufficiently large $t$, put $L:=\log t$ and prescribe
		\begin{equation}\label{eq:non_L1_prescribe}
			x(t):=\frac{1+L^{-1/2}}{L}, \qquad g(t):=x'(t)+f(x(t)).
		\end{equation}
		Then $x(t)/u(t)\to 1$ as $t\to\infty$. Differentiating gives
		\[
		x'(t)=-\frac{1}{tL^2}\left(1+\frac{3}{2\sqrt{L}}\right),
		\]
		whereas
		\[
		f(x(t))=\frac{(1+L^{-1/2})^2}{L^2}\exp\!\left(-\frac{L}{1+L^{-1/2}}\right).
		\]
		Hence
		\[
		\frac{f(x(t))}{|x'(t)|}
		= \frac{(1+L^{-1/2})^2}{1+3/(2\sqrt{L})}\exp\!\left(\frac{\sqrt{L}}{1+L^{-1/2}}\right)
		\longrightarrow \infty,
		\]
		which implies $g(t)>0$ for all sufficiently large $t$. Moreover,
		$f(x(t))\sim e^{-1} e^{\sqrt{L}}/(tL^2)$ as $t\to\infty$, so the substitution
		$L=\log t$ gives $\int^\infty f(x(t))\,dt = \infty$. Since $x(t)\to 0$ and
		$x'(t)<0$, integrating \eqref{eq:non_L1_prescribe} over $[t,R]$ gives
		\[
		\int_t^R g(s)\,ds = x(R)-x(t)+\int_t^R f(x(s))\,ds \longrightarrow \infty
		\qquad \text{as } R\to\infty,
		\]
		so $\int_t^\infty g(s)\,ds=\infty$ for every sufficiently large $t$.
	\end{proof}

\end{document}